\documentclass[twoside,11pt]{article}

\usepackage[abbrvbib, preprint]{jmlr2e}

\newcommand{\norm}[1]{\left|\left| #1 \right|\right|}

\renewcommand{\hat}{\widehat}
\renewcommand{\tilde}{\widetilde}

\allowdisplaybreaks

\newtheorem{remark*}{Remark}
\newtheorem{innercustomthm}{Theorem}
\newenvironment{customthm}[1]
{\renewcommand\theinnercustomthm{#1}\innercustomthm}
{\endinnercustomthm}

\newtheorem{innercustomprop}{Proposition}
\newenvironment{customprop}[1]
{\renewcommand\theinnercustomprop{#1}\innercustomprop}
{\endinnercustomprop}

\jmlrheading{1}{2021}{1-25}{4/00}{10/00}{zhang21}{Yikun Zhang and Yen-Chi Chen}

\ShortHeadings{The EM Perspective of DMS Algorithm}{Zhang and Chen}
\firstpageno{1}

\begin{document}
	
	\title{The EM Perspective of Directional Mean Shift Algorithm}
	
	\author{\name Yikun Zhang \email yikun@uw.edu
		\AND
		\name Yen-Chi Chen \email yenchic@uw.edu \\
		\addr University of Washington, Seattle, WA 98195}
	
	\editor{}
	
	\maketitle
	
\begin{abstract}%
	The directional mean shift (DMS) algorithm is a nonparametric method for pursuing local modes of densities defined by kernel density estimators on the unit hypersphere. In this paper, we show that any DMS iteration can be viewed as a generalized Expectation-Maximization (EM) algorithm; in particular, when the von Mises kernel is applied, it becomes an exact EM algorithm. Under the (generalized) EM framework, we provide a new proof for the ascending property of density estimates and demonstrate the global convergence of directional mean shift sequences. Finally, we give a new insight into the linear convergence of the DMS algorithm.
\end{abstract}


\section{Introduction}

The directional mean shift (DMS) algorithm is a generalization of the regular mean shift algorithm with Euclidean data \citep{MS1975,MS1995,MS2002,carreira2015review} toward handling directional data, which are assumed to consist of observations lying on a unit hypersphere $\Omega_q= \{\bm{x}\in\mathbb{R}^{q+1}: \norm{\bm{x}}_2 = 1\}$. Analogous to the regular mean shift procedure, the DMS algorithm is built upon a directional kernel density estimator (KDE) \citep{KDE_Sphe1987,KDE_direct1988,Zhao2001,Exact_Risk_bw2013} and conducts fixed-point iterations to pursue the local modes of the directional KDE.

Several different forms of the DMS algorithm have been proposed and applied to various fields such as gene clustering \citep{Multi_Clu_Gene2005}, medical structure classification \citep{DMS_topology2010}, seismological analysis \citep{DirMS2020} in the last two decades. Compared to its success in practice, the convergence theory of the DMS algorithm is less developed. The convergence of density estimates along any mean shift iterative sequence (or the so-called ascending property) was proved by \cite{MS1995,MS2002,MS2007_pf} for the regular mean shift algorithm and generalized to the mean shift algorithm on arbitrary manifolds by \cite{Nonlinear_MS_man2009}. This ascending property was specialized to a DMS algorithm with the von Mises kernel by \cite{vMF_MS2010}. \cite{DMS_topology2010} and \cite{DirMS2020} later showed that the DMS algorithm's ascending property holds for any convex and differentiable directional kernel. 

However, the convergence of DMS sequences is a more challenging problem. Such convergence has been well-studied in the regular mean shift scenarios \citep{MS1995,MS2007_pf,MS_EM2007,MS_onedim2013,MS2015_Gaussian,Ery2016}. In the directional data setting, the only work that we are aware of is \cite{DirMS2020}, but their result is limited to small neighborhoods of local modes. The global behavior of DMS sequences remains unclear.  

To address the global convergence of any DMS sequence on $\Omega_q$, we formulate the DMS algorithm into a special case of the famous (generalized) Expectation-Maximization (EM) algorithm \citep{EM1977,EM_Jeff1983,EM_Extensions2008,EM2017}. It allows us to present a new proof of the ascending property and study the convergence of DMS sequences based on extensive research in the (generalized) EM algorithm. We are inspired by \cite{MS_EM2007} for his method of writing the regular mean shift with Euclidean data as a (generalized) EM algorithm. 
However, the formulation is more difficult in any DMS algorithm because directional data lie on a nonlinear manifold. 
We overcome the problem by projecting the data from the original hypersphere $\Omega_q$ onto a high dimensional hypersphere $\Omega_{q+1}$
and design a mixture model to associate the DMS algorithm with the (generalized) EM algorithm. \\

\noindent\textbf{Main results} $\,$ The main objective of this paper is to show that the usual DMS algorithm (reviewed in Section~\ref{Sec:DMS}) is a (generalized) EM algorithm. Specially, our work contributes the following:
\begin{enumerate}
	\item We construct a maximum likelihood problem viewing the directional KDE as a mixture model with a single parameter $\bm{\mu}\in \Omega_q$ in a higher dimensional sphere $\Omega_{q+1}$. Fitting the mixture model with a (generalized) EM algorithm at a single point $\bm{Y}_1=(0,...,0,1)^T\in \Omega_{q+1}$ leads to the DMS algorithm (Section~\ref{Sec:DMS_EM}).
	\item Under the (generalized) EM framework, we prove the ascending property of the DMS algorithm based on the monotonicity of observed likelihoods along any (generalized) EM iteration (Section~\ref{Sec:ascending}).
	\item We derive the global convergence of DMS sequences using properties from the EM perspective (Section~\ref{Sec:conv}).
	\item Finally, compared to the result in \cite{DirMS2020}, we provide an alternative but straightforward proof of linear convergence of the DMS algorithm (Section~\ref{Sec:RC_DMS}). 
\end{enumerate}

\noindent \textbf{Notation} $\,$ Let $\Omega_q=\left\{\bm{x}\in \mathbb{R}^{q+1}: \norm{\bm{x}}_2=1 \right\}$ be the $q$-dimensional unit sphere in the ambient Euclidean space $\mathbb{R}^{q+1}$, where $\norm{\cdot}_2$ is the usual Euclidean norm (or $L_2$-norm) in $\mathbb{R}^{q+1}$. The inner product between two vectors $\bm{x},\bm{y} \in \mathbb{R}^{q+1}$ is denoted by $\bm{x}^T\bm{y}$. Given a smooth function $f: \mathbb{R}^{q+1} \to \mathbb{R}$, we denote its total gradient and Hessian by $\nabla f(\bm{x}) = \left(\frac{\partial f(\bm{x})}{\partial x_1},..., \frac{\partial f(\bm{x})}{\partial x_{q+1}} \right)^T \in \mathbb{R}^{q+1}$ and 
\[
\nabla\nabla f(\bm{x}) = \begin{pmatrix}
\frac{\partial^2 f(\bm{x})}{\partial x_1^2} & \cdots &\frac{\partial^2 f(\bm{x})}{\partial x_1 \partial x_{q+1}}\\
\vdots & \ddots & \vdots\\
\frac{\partial^2 f(\bm{x})}{\partial x_{q+1}\partial x_1} & \cdots & \frac{\partial^2 f(\bm{x})}{\partial x_{q+1}^2}\\
\end{pmatrix} \in \mathbb{R}^{(q+1)\times (q+1)},
\]
respectively. We use the big-$\mathcal{O}$ notation $p_t = \mathcal{O}(q_t)$ for a sequence of vectors $p_t$ and positive scalars $q_t$ if there is a global constant $C$ such that $\norm{p_t}_2 \leq C q_t$ for all sufficiently large $t$.

\section{Background}
\label{Sec:BG}

\subsection{Kernel Density Estimation with Directional Data}

Let $\bm{X}_1,...,\bm{X}_n \in \Omega_q\subset\mathbb{R}^{q+1}$ be a random sample generated from a directional density function $f$ on $\Omega_q$ with $\int_{\Omega_q} f(\bm{x}) \, \omega_q(d\bm{x})=1,$ 
where $\omega_q$ is the Lebesgue measure on $\Omega_q$. The directional KDE at point $\bm{x}\in \Omega_q$ is often written as \citep{Beran1979,KDE_Sphe1987,KDE_direct1988,Exact_Risk_bw2013}:
\begin{equation}
\label{Dir_KDE}
\hat{f}_h(\bm{x}) = \frac{c_{h,q,L}}{n} \sum_{i=1}^n L\left(\frac{1-\bm{x}^T \bm{X}_i}{h^2} \right),
\end{equation}
where $L$ is a directional kernel (i.e., a rapidly decaying function with nonnegative values and defined on $[0,\infty) \subset \mathbb{R}$), $h>0$ is the bandwidth parameter, and $c_{h,q,L}$ is a normalizing constant satisfying
\begin{equation}
\label{asym_norm_const}
c_{h,q,L}^{-1} = \int_{\Omega_q} L\left(\frac{1-\bm{x}^T \bm{y}}{h^2} \right) \omega_q(d\bm{y}) = \mathcal{O}(h^q).
\end{equation}

The bandwidth selection is crucial in determining the performance of directional KDEs because it controls the bias-variance tradeoff. 
There are various reliable bandwidth selection mechanisms in the literature \citep{KDE_Sphe1987,KDE_direct1988,Auto_bw_cir2008,KDE_torus2011,Oliveira2012,Exact_Risk_bw2013,Nonp_Dir_HDR2020}. 
Compared to the bandwidth, the choice of the kernel is less influential. One popular candidate is the so-called von Mises kernel $L(r) = e^{-r}$, whose name originates from the famous $q$-von Mises-Fisher (vMF) distribution on $\Omega_q$ with the following density: 
\begin{equation}
\label{vMF_density}
f_{\text{vMF}}(\bm{x}|\bm{\mu},\kappa) = C_q(\kappa) \cdot \exp\left(\kappa \bm{\mu}^T \bm{x} \right) \quad \text{ and } \quad C_q(\kappa) = \frac{\kappa^{\frac{q-1}{2}}}{(2\pi)^{\frac{q+1}{2}} \mathcal{I}_{\frac{q-1}{2}}(\kappa)},
\end{equation}
where $\bm{\mu} \in \Omega_q$ is the directional mean, $\kappa \geq 0$ is the concentration parameter, and $\mathcal{I}_{\alpha}(\kappa)$ is the modified Bessel function of the first kind of order $\kappa$; see \cite{Stat_Dir_Data1975,Mardia2000directional} for more details. With the von Mises kernel, the directional KDE \eqref{Dir_KDE} becomes a mixture of $q$-von Mises-Fisher densities as:
\begin{equation}
\label{Dir_KDE_vMF}
\hat{f}_h(\bm{x}) = \frac{1}{n} \sum_{i=1}^n f_{\text{vMF}}\left(\bm{x};\bm{X}_i,\frac{1}{h^2} \right) = \frac{C_q\left(\frac{1}{h^2} \right)}{n} \sum_{i=1}^n \exp\left(\frac{\bm{x}^T \bm{X}_i}{h^2} \right).
\end{equation}
The EM framework for estimating the parameters of a mixture of vMF distributions is also well-studied. We offer a detailed review in Appendix \ref{Appendix:Detail_EM_vMF}. Another potential choice is the following truncated convex kernel proposed by \cite{MSBC_Cir2012,MSC_Dir2014}:
\begin{equation}
\label{convex_kernel}
L(r) = 
\begin{cases}
(1-r)^p & \text{ if } 0\leq r \leq 1,\\
0   & \text{ otherwise,}
\end{cases}
\end{equation}
for some integer $p\geq 1$.

\subsection{Mean Shift Algorithm with Directional Data}
\label{Sec:DMS}

We follow the technique in \cite{MSC_Dir2014} to derive the most commonly used directional mean shift (DMS) algorithm. Assume that the kernel $L$ is differentiable except for finitely many points on $[0,\infty)$. Given the directional KDE $\hat{f}_h(\bm{x})$ in \eqref{Dir_KDE}, we introduce a Lagrangian multiplier $\lambda$ to maximize $\hat{f}_h(\bm{x})$ under the constraint $\bm{x}^T\bm{x}=1$ as follows:
$$\hat{\mathcal{L}}_h(\bm{x},\lambda)= \frac{c_{h,q,L}}{n} \sum_{i=1}^n L\left(\frac{1-\bm{x}^T \bm{X}_i}{h^2} \right) + \lambda(1-\bm{x}^T\bm{x}).$$
Taking the partial derivatives of $\hat{\mathcal{L}}_h$ with respect to $\bm{x}$ and $\lambda$ and setting them to zero yield that
$$\frac{\partial \hat{\mathcal{L}}_h}{\partial \bm{x}} = -\frac{c_{h,q,L}}{nh^2} \sum_{i=1}^n \bm{X}_iL'\left(\frac{1-\bm{x}^T \bm{X}_i}{h^2} \right) -2\lambda\bm{x}=0 \quad \text{ and } \quad 1-\bm{x}^T\bm{x}=0.$$
Solving for $\bm{x}$ leads to the fixed-point iteration equation $\bm{x}^{(t+1)} = F(\bm{x}^{(t)})$ for the DMS algorithm, where
\begin{equation}
\label{Dir_MS}
F(\bm{x}) = -\frac{\sum_{i=1}^n \bm{X}_iL'\left(\frac{1-\bm{x}^T\bm{X}_i}{h^2} \right)}{\norm{\sum_{i=1}^n \bm{X}_iL'\left(\frac{1-\bm{x}^T\bm{X}_i}{h^2} \right)}_2}.
\end{equation}
\cite{DirMS2020} suggested an alternative derivation of the DMS algorithm based on the following equality on $\Omega_q$:
$$\hat{f}_h(\bm{x}) = \tilde{f}_h(\bm{x}) \quad \text{ with } \quad \tilde{f}_h(\bm{x}) = \frac{c_{h,q,L}}{n} \sum_{i=1}^n L\left(\frac{1}{2}\norm{\frac{\bm{x}-\bm{X}_i}{h}}_2^2 \right),$$
and wrote out an explicit formula for the mean shift vector. Based on their derivation, the preceding fixed-point iteration scheme \eqref{Dir_MS} can be written in terms of $\nabla\hat{f}_h(\bm{x})$ as:
\begin{equation}
\label{Dir_MS_fix}
\bm{x}^{(t+1)} = F(\bm{x}^{(t)})= \frac{\nabla\hat{f}_h(\bm{x}^{(t)})}{\norm{\nabla\hat{f}_h(\bm{x}^{(t)})}_2},
\end{equation}
where $\nabla$ is the total gradient operator in the ambient space $\mathbb{R}^{q+1}$. 

\section{Directional Mean Shift as a Generalized EM Algorithm}
\label{Sec:DMS_EM}

\subsection{Detailed Derivations}
\label{Sec:Detailed_deri}

Recall that $\{\bm{X}_1,...,\bm{X}_n\} \subset \Omega_q$ is the observed directional dataset by the DMS algorithm in Section~\ref{Sec:DMS}. 
Our derivation starts with defining a mixture model on \emph{a higher dimensional hypersphere $\Omega_{q+1}$}. 
Consider the following mixture model for any point $\bm{y}\in\Omega_{q+1}$:
\begin{equation}
\label{Mixture_Density}
f_{L,q+1}(\bm{y}) = \sum_{i=1}^n \alpha_i \cdot C_{\kappa_i,q+1,L} \cdot L\left[\kappa_i(1-\bm{y}^T \bm{\nu}_i) \right],
\end{equation}
where 
$\{\alpha_i,\bm{\nu}_i,\kappa_i\}_{i=1}^n$ are the parameters of this mixture model
such that
$\sum\limits_{i=1}^n \alpha_i=1$, $\alpha_i \geq 0$ for $i=1,...,n$, and $C_{\kappa_i,q+1,L}$ is the normalizing constant such that 
\begin{equation}
\label{norm_const}
C_{\kappa_i,q+1,L}^{-1} = \int_{\Omega_{q+1}} L\left[\kappa_i\left(1-\bm{y}^T\bm{\nu}_i \right)\right] \omega_{q+1}(d\bm{y})= w_q(\Omega_q) \int_{-1}^1 L\left[\kappa_i(1-t)\right] \cdot (1-t^2)^{\frac{q}{2}-1} dt
\end{equation}
with $\omega_q(\Omega_q)$ as the surface area of $\Omega_q$. 
Equation \eqref{norm_const} 
demonstrates that the normalizing constant $C_{\kappa_i,q+1,L}$ is independent of the choice of the directional mean parameter $\bm{\nu}_i$. The parameters $\alpha_i$  and $\kappa_i$ represent the weight (mixing proportion) 
and the amount of smoothing (concentration) applied to the $i$-th observation, respectively.
Later, we will take $\alpha_i = \frac{1}{n}$
and $\kappa_i = \frac{1}{h^2}$, though they would depend on the index $i$
in a generic case.
One should keep in mind that  $\alpha_i,\kappa_i$ are fixed and given.  

We introduce an unknown parameter $\bm{\mu}$ relating
the directional mean parameter  $\bm{\nu}_i$ in the mixture model \eqref{Mixture_Density} to each observation $\bm{X}_i$ in the directional dataset as:
\begin{equation}
\label{mean_pars}
\bm{\nu}_i = \bm{\nu}_i(\bm{\mu}) = \left(\bm{X}_i\sqrt{1-\left(\bm{\mu}^T\bm{X}_i \right)^2}, \bm{\mu}^T\bm{X}_i \right)^T \in \Omega_{q+1}.
\end{equation}
Given $\bm{\mu}\in \Omega_q$, the distribution in \eqref{Mixture_Density} is a mixture of $n$ densities on $\Omega_{q+1}$
with 
directional means $\{\bm{\nu}_i(\bm{\mu})\}_{i=1}^n$, 
\emph{fixed} mixing proportions $\{\alpha_i\}_{i=1}^n$
and \emph{fixed} concentration parameters $\{\kappa_i\}_{i=1}^n$. Varying the parameter $\bm{\mu} \in \Omega_{q+1}$ in the mixture model leads to a spherical transformation of the whole mixture distribution in a higher dimensional unit sphere $\Omega_{q+1}$. 
Our goal is to find the maximum likelihood estimate (MLE) of $\bm{\mu}$ given a single observation from the mixture model.

We consider the following hypothetical setup for deriving the EM algorithm and connecting it with the DMS algorithm.
Suppose that we are given $\bm{Y}_1$, an observation in $\Omega_{q+1}$. 
We fit the mixture model \eqref{Mixture_Density} to $\bm{Y}_1$
and attempt to find the MLE of $\bm{\mu}$ based on this single observation.
For such mixture model, 
finding the MLE of $\bm{\mu}$
can be done via the EM algorithm. 
To this end, we introduce a hidden/latent random variable $Z_1$
indicating from which component $\bm{Y}_1$ is generated. Namely, $Z_1=i$ if $\bm{Y}_1$ is sampled from $C_{\kappa_i,q+1,L} \cdot L\left[\kappa_i(1-\bm{y}^T \bm{\nu}_i) \right]$ in \eqref{Mixture_Density}.
The EM framework for obtaining the MLE of $\bm{\mu}$ given the observed data $\bm{Y}_1$ is as follows:

\noindent $\bullet$ {\bf E-Step}: Assuming the values of $Z_1$ is known, the complete log-likelihood is written as:
\begin{equation}
\label{complete_ll_MS}
\log P(\bm{Y}_1,{Z}_1|\bm{\mu}) = \sum_{i=1}^n \mathbbm{1}_{\{Z_1=i\}} \cdot \left\{\log \alpha_i + \log C_{\kappa_i,q+1,L} + \log L\left[\kappa_i(1-\bm{Y}_1^T\bm{\nu}_i)\right] \right\}.
\end{equation}
Given $(\bm{Y}_1,\bm{\mu}^{(t)})$ at the $t$-th iteration, we take the expectation of the complete log-likelihood \eqref{complete_ll_MS} with respect to the current posterior distribution of ${Z}_1|\left(\bm{Y}_1,\bm{\mu}^{(t)} \right)$ as:
\begin{align}
\label{Q_function_MS}
\begin{split}
Q(\bm{\mu}|\bm{\mu}^{(t)}) &= \mathbb{E}_{{Z}_1|(\bm{Y}_1,\bm{\mu}^{(t)})}\left[\log P(\bm{Y}_1,{Z}_1|\bm{\mu}) \right]\\
&= \sum_{i=1}^n \left[\log \alpha_i + \log C_{\kappa_i,q+1,L} \right] \cdot P(Z_1=i|\bm{Y}_1,\bm{\mu}^{(t)}) \\
&\quad +\sum_{i=1}^n \log L\left[\kappa_i (1-\bm{Y}_1^T\bm{\nu}_i) \right] \cdot P(Z_1=i|\bm{Y}_1,\bm{\mu}^{(t)}).
\end{split}
\end{align}
By Bayes' theorem, the posterior distribution $P(Z_1=i|\bm{Y}_1,\bm{\mu}^{(t)})$ is computed as
\begin{align}
\label{post_hidden_MS}
\begin{split}
P(Z_1=i|\bm{Y}_1,\bm{\mu}^{(t)}) &= \frac{P(\bm{Y}_1|Z_1=i,\bm{\mu}^{(t)}) \cdot P(Z_1=i|\bm{\mu}^{(t)})}{P(\bm{Y}_1|\bm{\mu}^{(t)})}\\
&= \frac{\alpha_i \cdot C_{\kappa_i,q+1,L} \cdot L\left[\kappa_i\left(1-\bm{Y}_1^T\bm{\nu}_i^{(t)} \right) \right]}{\sum_{\ell=1}^n \alpha_{\ell} \cdot C_{\kappa_{\ell},q+1,L}\cdot L\left[\kappa_{\ell} \left(1-\bm{Y}_1^T \bm{\nu}_{\ell}^{(t)} \right) \right]}.
\end{split}
\end{align}

\noindent $\bullet$ {\bf M-Step}: 
To derive the M-step, we focus on a particular case where the observation $\bm{Y}_1=(0,...,0,1)^T \in \Omega_{q+1}$.
Since $\bm{\nu}_i^T \bm{Y}_1 = \bm{\mu}^T\bm{X}_i$ for $i=1,...,n$, the Q-function \eqref{Q_function_MS} reduces to
\begin{align}
\label{Q_function_MS2}
\begin{split}
Q(\bm{\mu}|\bm{\mu}^{(t)}) &= \sum_{i=1}^n \left[\log \alpha_i + \log C_{\kappa_i,q+1,L} \right] \cdot P(Z_1=i|\bm{Y}_1,\bm{\mu}^{(t)})\\
&\quad + \sum_{i=1}^n \log L\left[\kappa_i(1-\bm{\mu}^T\bm{X}_i) \right] \cdot P(Z_1=i|\bm{Y}_1,\bm{\mu}^{(t)}).
\end{split}
\end{align}
As is shown in \eqref{norm_const} that each normalizing constant $C_{\kappa_i,q+1,L}$ is independent of the directional mean parameter $\bm{\nu}_i$ for $i=1,...,n$ (and consequently, the parameter $\bm{\mu}$), we only need to maximize the second term in \eqref{Q_function_MS2} with respect to $\bm{\mu}$ under the constraint $\bm{\mu}^T\bm{\mu}=1$. To this end, we introduce a Lagrangian multiplier $\lambda$ and obtain the following Lagrangian function:
$$\mathcal{L}(\bm{\mu},\lambda) = \sum_{i=1}^n \log L\left[\kappa_i(1-\bm{\mu}^T\bm{X}_i) \right] \cdot P(Z_1=i|\bm{Y}_1,\bm{\mu}^{(t)}) + \lambda \left(1-\bm{\mu}^T\bm{\mu} \right).$$
Taking the partial derivatives of $\mathcal{L}$ with respect to $\bm{\mu}$ and $\lambda$ and setting them to zeros yield that
$$\frac{\partial \mathcal{L}}{\partial \bm{\mu}}= \sum_{i=1}^n P(Z_1=i|\bm{Y}_1,\bm{\mu}^{(t)}) \cdot \frac{-\kappa_i \bm{X}_i \cdot L'\left[\kappa_i(1-\bm{\mu}^T\bm{X}_i)\right]}{L\left[\kappa_i(1-\bm{\mu}^T\bm{X}_i)\right]} - 2\lambda \bm{\mu}=0 \quad \text{ and }\quad 1-\bm{\mu}^T\bm{\mu}=0.$$
This further implies that
\begin{align}
\label{mu_iter_MS}
\begin{split}
\hat{\bm{\mu}} &= \frac{\sum_{i=1}^n P(Z_1=i|\bm{Y}_1,\bm{\mu}^{(t)}) \cdot \frac{-\kappa_i \bm{X}_i \cdot L'\left[\kappa_i(1-\hat{\bm{\mu}}^T\bm{X}_i)\right]}{L\left[\kappa_i(1-\hat{\bm{\mu}}^T\bm{X}_i)\right]}}{\norm{\sum_{i=1}^n P(Z_1=i|\bm{Y}_1,\bm{\mu}^{(t)}) \cdot \frac{-\kappa_i \bm{X}_i \cdot L'\left[\kappa_i(1-\hat{\bm{\mu}}^T\bm{X}_i)\right]}{L\left[\kappa_i(1-\hat{\bm{\mu}}^T\bm{X}_i)\right]}}_2}\\
&= \frac{\sum_{i=1}^n \frac{\alpha_i \cdot C_{\kappa_i,q+1,L} \cdot L\left[\kappa_i\left(1-\bm{X}_i^T\bm{\mu}^{(t)} \right) \right]}{\sum_{\ell=1}^n \alpha_{\ell} \cdot C_{\kappa_{\ell},q+1,L}\cdot L\left[\kappa_{\ell} \left(1-\bm{X}_{\ell}^T \bm{\mu}^{(t)} \right) \right]} \cdot \frac{-\kappa_i \bm{X}_i \cdot L'\left[\kappa_i(1-\hat{\bm{\mu}}^T\bm{X}_i)\right]}{L\left[\kappa_i(1-\hat{\bm{\mu}}^T\bm{X}_i)\right]}}{\norm{\sum_{i=1}^n \frac{\alpha_i \cdot C_{\kappa_i,q+1,L} \cdot L\left[\kappa_i\left(1-\bm{X}_i^T\bm{\mu}^{(t)} \right) \right]}{\sum_{\ell=1}^n \alpha_{\ell} \cdot C_{\kappa_{\ell},q+1,L}\cdot L\left[\kappa_{\ell} \left(1-\bm{X}_{\ell}^T \bm{\mu}^{(t)} \right) \right]} \cdot \frac{-\kappa_i \bm{X}_i \cdot L'\left[\kappa_i(1-\hat{\bm{\mu}}^T\bm{X}_i)\right]}{L\left[\kappa_i(1-\hat{\bm{\mu}}^T\bm{X}_i)\right]}}_2}\\
&= \frac{-\sum_{i=1}^n \kappa_i \alpha_i \cdot C_{\kappa_i,q+1,L} \bm{X}_i \cdot L'\left[\kappa_i \left(1-\hat{\bm{\mu}}^T \bm{X}_i \right) \right] \cdot \frac{L\left[\kappa_i \left(1-\bm{X}_i^T\bm{\mu}^{(t)} \right) \right]}{L\left[\kappa_i\left( 1-\hat{\bm{\mu}}^T\bm{X}_i\right) \right]}}{\norm{\sum_{i=1}^n \kappa_i \alpha_i \cdot C_{\kappa_i,q+1,L} \bm{X}_i \cdot L'\left[\kappa_i \left(1-\hat{\bm{\mu}}^T \bm{X}_i \right) \right] \cdot \frac{L\left[\kappa_i \left(1-\bm{X}_i^T\bm{\mu}^{(t)} \right) \right]}{L\left[\kappa_i\left( 1-\hat{\bm{\mu}}^T\bm{X}_i\right) \right]}}_2},
\end{split}
\end{align}
where we plug in the posterior distribution \eqref{post_hidden_MS} in the second equality, and the term 
$$\sum_{\ell=1}^n \alpha_{\ell} \cdot C_{\kappa_{\ell},q+1,L}\cdot L\left[\kappa_{\ell} \left(1-\bm{X}_{\ell}^T \bm{\mu}^{(t)} \right) \right]$$
is independent of the outer summation index $i$ in both the numerator and denominator, so they are canceled out in the third equality. Therefore, equation \eqref{mu_iter_MS} defines a fixed-point iteration scheme for obtaining $\bm{\mu}^{(t+1)}$ in the EM framework, and an exact M-step will iterate \eqref{mu_iter_MS} until convergence.\\

Suppose that we take $\kappa_i=\frac{1}{h^2}$ and $\alpha_i=\frac{1}{n}$ for all $i=1,...,n$ in the fixed-point iteration equation \eqref{mu_iter_MS}. Then, the normalizing constant $C_{\kappa_i,q+1,L}$ will no longer depend on the summation index $i$, and the fixed-point iteration equation reduces to
\begin{equation}
\label{mu_fixed_point}
\hat{\bm{\mu}} = -\frac{\sum_{i=1}^n \bm{X}_i L'\left(\frac{1-\hat{\bm{\mu}}^T \bm{X}_i}{h^2} \right) \cdot \frac{L\left(\frac{1-\bm{X}_i^T \bm{\mu}^{(t)}}{h^2} \right)}{L\left(\frac{1-\hat{\bm{\mu}}^T \bm{X}_i}{h^2} \right)}}{\norm{\sum_{i=1}^n \bm{X}_i L'\left(\frac{1-\hat{\bm{\mu}}^T \bm{X}_i}{h^2} \right) \cdot \frac{L\left(\frac{1-\bm{X}_i^T \bm{\mu}^{(t)}}{h^2} \right)}{L\left(\frac{1-\hat{\bm{\mu}}^T \bm{X}_i}{h^2} \right)}}_2}
\end{equation}
or more specifically, 
\begin{equation}
\label{mu_fixed_point2}
\hat{\bm{\mu}}^{(k+1;t)} = -\frac{\sum_{i=1}^n \bm{X}_i L'\left(\frac{1- \bm{X}_i^T \hat{\bm{\mu}}^{(k;t)}}{h^2} \right) \cdot \frac{L\left(\frac{1-\bm{X}_i^T \bm{\mu}^{(t)}}{h^2} \right)}{L\left(\frac{1- \bm{X}^T_i \hat{\bm{\mu}}^{(k;t)}}{h^2} \right)}}{\norm{\sum_{i=1}^n \bm{X}_i L'\left(\frac{1-\bm{X}_i^T \hat{\bm{\mu}}^{(k;t)} }{h^2} \right) \cdot \frac{L\left(\frac{1-\bm{X}_i^T \bm{\mu}^{(t)}}{h^2} \right)}{L\left(\frac{1- \bm{X}_i^T \hat{\bm{\mu}}^{(k;t)}}{h^2} \right)}}_2}
\end{equation}
for $k=0,1,2,...$. Therefore, the complete EM algorithm consists of two nested loops: the outer one is the usual iteration over $t$, while the inner one iterates $\{\hat{\bm{\mu}}^{(k;t)}\}_{k=0}^{\infty}$ under the fixed $\bm{\mu}^{(t)}$ for an exact M-step. If we choose the initial value for iterating \eqref{mu_fixed_point2}  as
$\hat{\bm{\mu}}^{(0;t)} = \bm{\mu}^{(t)}$ and do a single iteration of \eqref{mu_fixed_point2}, i.e., $\bm{\mu}^{(t+1)} = \hat{\bm{\mu}}^{(1;t)}$,  we obtain
\begin{equation}
\label{EM_MS_iter}
\bm{\mu}^{(t+1)} = -\frac{\sum_{i=1}^n \bm{X}_i L'\left(\frac{1-\bm{X}_i^T\bm{\mu}^{(t)}}{h^2}\right)}{\norm{\sum_{i=1}^n \bm{X}_i L'\left(\frac{1-\bm{X}_i^T\bm{\mu}^{(t)}}{h^2}\right)}_2},
\end{equation}
which coincides with the mean shift updating formula with directional data (c.f., equation \eqref{Dir_MS}). 
In other words, the DMS algorithm can be viewed as the above EM algorithm with an inexact M-step update. 
More importantly, as we will demonstrate in Section \ref{Sec:ascending}, the Q-function $Q(\bm{\mu}|\bm{\mu}^{(t)})$ is non-decreasing along the iterative path defined by \eqref{EM_MS_iter} under some mild conditions on kernel $L$. Therefore, the DMS algorithm with a general kernel is indeed a generalized EM algorithm.

\begin{remark*}
	Our derivation is different from
	the suggestion given in \cite{Nonlinear_MS_man2009}. 
	In our derivation,
	we cannot  define a group action of $\mathcal{G}$ on $\Omega_q$ such that the componentwise density $C_{\kappa_i,q+1,L} \cdot L\left[\kappa_i(1-\bm{y}^T \bm{\nu}_i) \right]$ in \eqref{Mixture_Density} is equivalent to the original density $C_{\kappa_i,q,L}\cdot  L\left[\kappa_i(1-(g\cdot\bm{y})^T\bm{X}_i) \right]$ under some $g\in \mathcal{G}$ because these two densities lie on hyperspheres with different dimensions.
\end{remark*}

\subsection{Directional Mean Shift with von Mises Kernel as an Exact EM Algorithm}
\label{Sec:EM_MS_vM_kernel}

As the Gaussian mean shift with Euclidean data can be written as an EM algorithm \citep{MS_EM2007}, we show in this subsection that the DMS algorithm with the von Mises kernel is also an EM algorithm. That is, the preceding M-step in Section \ref{Sec:Detailed_deri} is exact when we apply the von Mises kernel. This is because the von Mises kernel $L(r)=e^{-r}$ satisfies the property that $L'\propto L$, and the above M-step \eqref{mu_iter_MS} is simplified as:
\begin{equation}
\label{mu_iter_MS_vMF}
\bm{\mu}^{(t+1)} = \frac{\sum_{i=1}^n \kappa_i\alpha_i C_{q+1}(\kappa_i) \bm{X}_i \cdot \exp\left(\kappa_i\bm{X}_i^T \bm{\mu}^{(t)} \right)}{\norm{\sum_{i=1}^n \kappa_i\alpha_i C_{q+1}(\kappa_i) \bm{X}_i \cdot \exp\left(\kappa_i\bm{X}_i^T \bm{\mu}^{(t)} \right)}_2},
\end{equation}
Now, if we take $\kappa_i=\frac{1}{h^2}$ and $\alpha_i=\frac{1}{n}$ for all $i=1,...,n$ as before, the exact M-step in \eqref{mu_iter_MS_vMF} becomes
\begin{equation}
\label{EM_MS_iter_vMF}
\bm{\mu}^{(t+1)} = \frac{\sum_{i=1}^n \bm{X}_i \exp\left(\frac{\bm{X}_i^T \bm{\mu}^{(t)}}{h^2} \right)}{\norm{\sum_{i=1}^n \bm{X}_i \exp\left(\frac{\bm{X}_i^T \bm{\mu}^{(t)}}{h^2} \right)}_2},
\end{equation}
which is precisely the DMS algorithm with the von Mises kernel. As a result, the DMS algorithm with the von Mises kernel is an exact EM algorithm.

\subsection{Comparison to the Regular Mean Shift Algorithm with Euclidean Data}

Our derivation is motivated by the discovery of \cite{MS_EM2007}, in which the author obtained an elegant connection between the regular mean shift algorithm with Euclidean data and (generalized) EM algorithm.
The key of the derivation in \cite{MS_EM2007} is to define an artificial maximum likelihood problem where
\begin{enumerate}
	\item the underlying distribution is a mixture model defined by the original KDE so that the means of its components are given by the dataset to which the mean shift algorithm is applied, but the mixture model contains only a single displacement parameter,
	\item and, the mixture model is fit on the sample with a single point, the origin $\bm{Y}=\bm{0}$.
\end{enumerate}

Under the above construction, \cite{MS_EM2007} proved that any maximum likelihood estimate of the displacement vector is a mode of the original KDE. However, as is argued in Appendix \ref{Appendix:EM_displacement}, naively applying the idea of \cite{MS_EM2007} would fail. We have to make a non-trivial extension of his approach to bridge the connection between DMS and EM algorithms.
The two key changes we made in our derivation are:
\begin{enumerate}
	\item 
	We construct a mixture model in a higher-dimensional hypersphere $\Omega_{q+1}$ and
	introduce the parameter $\bm{\mu}\in \Omega_q$ controlling the projections of observations $\{\bm{X}_1,...,\bm{X}_n\} \subset \Omega_q$
	onto $\Omega_{q+1}$.
	\item We fit the mixture model on the sample containing a single point $\bm{Y}_1=(0,...,0,1)^T \in \Omega_{q+1}$.
\end{enumerate}

\section{Convergence of Directional Mean Shift as a Generalized EM Algorithm}
\label{Sec:Conv_DMS_as_EM}

In this section, we show that the Q-function \eqref{Q_function_MS2} is non-decreasing along any DMS sequence \eqref{EM_MS_iter} under some regularity conditions on kernel $L$. As the directional KDE will be identical to the corresponding observed likelihood up to a multiplicative constant, we furnish a new proof of the ascending property of density estimates along DMS sequences based on the generalized EM perspective, which is different from the existing result (e.g., \cite{vMF_MS2010,DMS_topology2010,DirMS2020}). Moreover, the connection to the (generalized) EM algorithm implies that the DMS sequence converges from almost any starting point on $\Omega_q$, though it may be stuck in a saddle point or a local minimum of the directional KDE $\hat{f}_h$ in \eqref{Dir_KDE}; see such non-typical behaviors of EM algorithms in Section 3.6.1 and 3.6.2 of \cite{EM_Extensions2008}. Nevertheless, the latter cases are unlikely in practice because both saddle points and minima are unstable for maximization \citep{MS_EM2007}. Therefore, the convergence of DMS sequences will almost always be to a local mode of $\hat{f}_h$ in practice. This conclusion is much stronger than the local convergence in Theorem 11 of \cite{DirMS2020}.

\subsection{Assumptions}
\label{Sec:assump}

To guarantee the convergence of the DMS iteration \eqref{EM_MS_iter} as a generalized EM algorithm, we make the following assumption on the directional KDE \eqref{Dir_KDE} and kernel $L$:
\begin{itemize}
	\item {\bf (C1)} The number of local modes of $\hat{f}_h$ on $\Omega_q$ is finite, and the modes are isolated from other critical points.
	\item {\bf (C2)} $L: [0,\infty) \to [0,\infty)$ is non-increasing, differentiable, and convex with $0<L(0) < \infty$.
\end{itemize}

Condition (C1) is commonly assumed for the convergence of regular mean shift sequences \citep{MS2007_pf,MS2015_Gaussian}. It is an evident corollary for the regular KDE with various kernels \citep{Silverman1981,hall2004,Non_mode_infer2016,Modal_Age_stat2020}, and one can easily extend the result to directional KDEs. Condition (C2) is a natural assumption on kernel $L$; see, for instance, \cite{MS2002,MS2007_pf,DMS_topology2010,DirMS2020}. It is satisfied by many directional kernels, such as the von Mises kernel $L(r)=e^{-r}$. The differentiability of kernel $L$ can be relaxed to include finitely many non-differentiable points on $[0,\infty)$. Such relaxation enables our convergence results to incorporate the DMS algorithms with kernels having bounded supports (e.g., kernels in \eqref{convex_kernel}).

\subsection{Ascending Property of Directional Mean Shift Algorithm}
\label{Sec:ascending}

The key reason why the DMS iteration is a generalized EM algorithm is that a single iteration of every incomplete M-step with the initial value $\bm{\mu}^{(t)}$ as in \eqref{EM_MS_iter} ensures the following non-decreasing property of the Q-function \eqref{Q_function_MS2}.

\begin{theorem}
	\label{Q_fun_increase}
	Under condition (C2), the Q-function \eqref{Q_function_MS2} satisfies $Q(\bm{\mu}^{(t+1)}|\bm{\mu}^{(t)}) \geq Q(\bm{\mu}^{(t)}|\bm{\mu}^{(t)})$ 
	where $\bm{\mu}^{(t+1)} = -\frac{\sum_{i=1}^n \bm{X}_i L'\left(\frac{1-\bm{X}_i^T\bm{\mu}^{(t)}}{h^2}\right)}{\norm{\sum_{i=1}^n \bm{X}_i L'\left(\frac{1-\bm{X}_i^T\bm{\mu}^{(t)}}{h^2}\right)}_2}$ is from the iteration \eqref{EM_MS_iter}.
\end{theorem}

The essences of the proof of Theorem~\ref{Q_fun_increase} are the following two inequalities:
\begin{equation}
\label{Q_fun_increase_pf}
\log y - \log x \geq \frac{y-x}{x} \quad \text{ and } \quad L(y)-L(x) \geq L'(x) \cdot (y-x),
\end{equation}
which are guaranteed by the concavity of $\log(\cdot)$ as well as the convexity and differentiability of the kernel $L$; see Appendix~\ref{Appendix:Q_fun_increase_pf}. The condition (C2) on kernel $L$ provides a sufficient condition under which the directional mean shift iteration does increase the Q-function \eqref{Q_function_MS2}, but it is by no means necessary. One can choose a different kernel and show that the iteration \eqref{EM_MS_iter} increases the Q-function using a similar argument but under other conditions (e.g., $\log L(r)$ is convex).

The observed log-likelihood with data point $\bm{Y}_1=(0,...,0,1)^T\in \Omega_{q+1}$ under the EM framework of Section~\ref{Sec:Detailed_deri} embraces the following association with the logarithm of the directional KDE $\hat{f}_h$ in \eqref{Dir_KDE}.

\begin{proposition}
	\label{obs_ll_Dir_KDE}
	The observed log-likelihood $\log P(\bm{Y}_1|\bm{\mu})$ is given by 
	\begin{equation}
	\label{obs_log_likelihood}
	\log P(\bm{Y}_1|\bm{\mu}) = \log\left[\sum_{\ell=1}^n \alpha_{\ell} \cdot C_{\kappa_{\ell},q+1,L} \cdot L\left[\kappa_{\ell} \left(1-\bm{X}_{\ell}^T\bm{\mu} \right)\right] \right].
	\end{equation}
	In particular, when $\alpha_{\ell}=\frac{1}{n}$ and $\kappa_{\ell}=\frac{1}{h^2}$ for all $\ell=1,...,n$ (i.e., simplifying to the DMS case), it implies that
	\begin{equation}
	\label{obs_ll_log_density}
	\log P(\bm{Y}_1|\bm{\mu}) = \log \hat{f}_h(\bm{\mu}) + \log \frac{C_{1/h^2,q+1,L}}{c_{h,q,L}},
	\end{equation}
	where $c_{h,q,L}$ is defined in \eqref{asym_norm_const} and $C_{\kappa_i,q+1,L}$ is elucidated in \eqref{norm_const}.
\end{proposition}

The proof of Proposition~\ref{obs_ll_Dir_KDE} is deferred to Appendix~\ref{Appendix:obs_ll_Dir_KDE}. Although some careful readers may note that $c_{h,q,L}= C_{1/h^2,q,L}$, the extra term $\log \frac{C_{1/h^2,q+1,L}}{c_{h,q,L}}$ is still irreducible in \eqref{obs_ll_log_density}, because we define the mixture density \eqref{Mixture_Density} for the (generalized) EM algorithm on a higher dimensional sphere $\Omega_{q+1}$. However, the presence of this extra term will not affect the ascending property of density estimates along DMS (or equivalently, generalized EM) sequences, as stated in the theorem below. 

\begin{theorem}
	\label{EM_Ascending}
	Let $\left\{\bm{\mu}^{(t)} \right\}_{t=0}^{\infty}$ be the path of successive points defined by the DMS algorithm \eqref{Dir_MS}. Under condition (C2), the sequence $\left\{\hat{f}_h(\bm{\mu}^{(t)}) \right\}_{t=0}^{\infty}$ is non-decreasing and thus converges.
\end{theorem}

Theorem~\ref{EM_Ascending} is the well-known monotonicity of observed (log-)likelihood sequences defined by the (generalized) EM algorithm (c.f., Theorem 1 in \cite{EM1977}). The proof is in Appendix~\ref{Appendix:EM_Ascending}.

\subsection{Convergence of Directional Mean Shift Sequences}
\label{Sec:conv}

The connection to the (generalized) EM algorithm allows us to derive a global convergence property of the DMS algorithm. 
\begin{theorem}
	\label{EM_conv_new}
	Let $\left\{\bm{\mu}^{(t)} \right\}_{t=0}^{\infty}$ be the path of successive points defined by the DMS algorithm \eqref{Dir_MS}. Under conditions (C1) and (C2), the sequence $\left\{\bm{\mu}^{(t)} \right\}_{t=0}^{\infty}$ converges to a local mode of $\hat f_h$
	from almost all starting point $\bm{\mu}^{(0)}\in \Omega_{q}$.
\end{theorem}

The high-level idea of the proof is as follows.
For the convergence of (generalized) EM sequences to stationary points (in most cases, local maxima), \cite{EM_Jeff1983} assumed the compactness of the parameter space and differentiability of the observed log-likelihood within the parameter space. In addition, the set of local modes of the observed log-likelihood is required to be discrete in order to guarantee the convergence of any generalized EM sequence to its single stationary point. 
We detail out these assumptions and the resulting convergence result for generalized EM (or equally, DMS) sequences in Appendix~\ref{Appendix:EM_conv} and show that conditions (C1) and (C2) imply these assumptions. 

In practice, the set of local minima and saddle points of $\hat{f}_h$ will have zero Lebesgue measure on $\Omega_q$, so the convergence of DMS sequences will almost always be to a local mode; see Appendix~\ref{Appendix:Basin_of_attr} for basins of attraction (or convergence domains) yielded by DMS algorithms. Therefore, our result (Theorem~\ref{EM_conv_new}) is stronger than the convergence theory of \cite{DirMS2020}, in which they only proved the convergence of DMS sequences around each local mode of $\hat{f}_h$.

\section{Rate of Convergence of Directional Mean Shift Algorithm}
\label{Sec:RC_DMS}

In this section, we study the rate of convergence of the DMS algorithm. We have shown in Sections~\ref{Sec:DMS_EM} and \ref{Sec:Conv_DMS_as_EM} that the DMS algorithm is a generalized EM algorithm and converges to local modes of the directional KDE from almost every starting point on $\Omega_q$. It is well-known that the convergence rate of an EM algorithm is generally linear (Page 102 in \cite{EM_Extensions2008}). \cite{MS_EM2007} studied in detail the rate of convergence for Gaussian mean shift in the homoscedastic and isotropic case. More recently, \cite{DirMS2020} wrote the DMS into a gradient ascent algorithm with an adaptive step size on $\Omega_q$ and argued that the DMS algorithm will converge linearly with a sufficiently small bandwidth around some neighborhoods of local modes. We will investigate the algorithmic convergence rate of the DMS algorithm in an alternative but straightforward approach. Consider the following assumptions: 
\begin{itemize}
	\item {\bf (D1)} The true directional density $f$ is continuously differentiable, and its partial derivatives are square integrable on $\Omega_q$.
	\item {\bf (D2)} Besides (C2), we assume that $L: [0,\infty) \to [0,\infty)$ satisfies 
	$$0< \int_0^{\infty} L(r)^k r^{\frac{q}{2}-1} dr < \infty \quad \text{ and } \quad 0< \int_0^{\infty} L'(r)^k r^{\frac{q}{2}-1} dr < \infty $$
	for all $q\geq 1$ and $k=1,2$, and
	$\int_{\epsilon^{-1}}^{\infty} \left|L'(r) \right| r^{\frac{q}{2}} dr = O(\epsilon)$
	as $\epsilon \to 0$ for any $\epsilon >0$.
\end{itemize}

Conditions (D1) and (D2) are common assumptions in directional data to establish the pointwise and uniform consistency of directional KDE $\hat{f}_h$ to the true density $f$ on $\Omega_q$ \citep{KDE_Sphe1987,KLEMELA2000,Zhao2001,Dir_Linear2013,Exact_Risk_bw2013,DirMS2020,alonso2020nonparametric}. 

Now, let $\bm{m}$ be a local mode of directional KDE $\hat{f}_h$ to which the DMS sequence $\{\bm{x}^{(t)}\}_{t=0}^{\infty}$ defined by \eqref{Dir_MS} converges. By Taylor's theorem and the fact that $\bm{m} =F(\bm{m})$, we deduce that
\begin{equation}
\label{DMS_fix_Taylor}
\bm{x}^{(t+1)} = F(\bm{x}^{(t)}) = \bm{m} + \nabla F(\bm{m}) \cdot (\bm{x}^{(t)} - \bm{m}) + \mathcal{O}\left(\norm{\bm{x}^{(t)} - \bm{m}}_2^2 \right),
\end{equation}
where $\nabla F(\bm{m}) \in \mathbb{R}^{(q+1)\times (q+1)}$ is the Jacobian of $F$ evaluated at $\bm{m}$. Thus, 
$\bm{x}^{(t+1)} - \bm{m} \approx \nabla F(\bm{m}) (\bm{x}^{(t)} - \bm{m})$
around a small neighborhood of $\bm{m}$ on $\Omega_q$. It implies that the rate of convergence is associated with the Jacobian $\nabla F(\bm{m})$ through its eigenvalue as follows:
\begin{equation}
\label{conv_rate_eigen}
\frac{\norm{\bm{x}^{(t+1)}-\bm{x}^*}_2}{\norm{\bm{x}^{(t)}-\bm{x}^*}_2} \leq \max_{i=1,...,q+1} |\Lambda_i|,
\end{equation}
where $\Lambda_1,...,\Lambda_{q+1}$ are the eigenvalues of $\nabla F(\bm{m})$. More importantly, the Jacobian $\nabla F(\bm{x})$ have the following useful properties.

\begin{proposition}
	\label{Jacobian_Prop}
	For any $\bm{x}\in \Omega_q \subset \mathbb{R}^{q+1}$,
	\begin{enumerate}[label=(\alph*)]
		\item the Jacobian $\nabla F(\bm{x})$ is a symmetric matrix in $\mathbb{R}^{(q+1)\times (q+1)}$, and has an explicit form
		\begin{equation}
		\label{Jacobian_F}
		\nabla F(\bm{x}) = \frac{\nabla\nabla \hat{f}_h(\bm{x})}{\norm{\nabla \hat{f}_h(\bm{x})}_2} - \frac{\nabla \hat{f}_h(\bm{x}) \nabla \hat{f}_h(\bm{x})^T \nabla\nabla \hat{f}_h(\bm{x})}{\norm{\nabla \hat{f}_h(\bm{x})}_2^3}.
		\end{equation}
		In particular, $\nabla F(\bm{m}) = \frac{(\bm{I}_{q+1}-\bm{m}\bm{m}^T) \nabla \nabla \hat{f}_h(\bm{m})}{\norm{\nabla \hat{f}_h(\bm{m})}_2}$ at any local mode $\bm{m}$ of $\hat{f}_h$.
		\item we have that $\max_{i=1,...,q+1}|\Lambda_i| \to 0$ as $h\to 0$ and $nh^q\to \infty$, where $\Lambda_1,...,\Lambda_{q+1}$ are the eigenvalues of $\nabla F(\bm{m})$.
	\end{enumerate}
\end{proposition}

The proof of Proposition~\ref{Jacobian_Prop} is in Appendix~\ref{Appendix:Jacobian_Prop_pf}. 
Note that $\norm{\nabla \hat{f}_h(\bm{x})}_2\to \infty$ as $h\to 0$ and $nh^q\to \infty$ for any $\bm{x}\in \Omega_q$; see Lemma 10 in \cite{DirMS2020}.
While $\norm{\nabla \hat{f}_h(\bm{x})}_2\to \infty$ may seem to be counter-intuitive, it is indeed a reasonable result because the total gradient $\nabla \hat{f}_h(\bm{x})$ in $\mathbb{R}^{q+1}$ consists of two parts: the gradient tangent to $\Omega_{q}$ (tangent gradient) and the gradient perpendicular to $\Omega_{q}$ (normal gradient). 
The tangent gradient is a stable quantity (since it is a well-defined quantity in terms of the true directional density $f$) and the normal gradient will diverge (since the true directional density $f$ does not have any normal gradient). 

Together with \eqref{conv_rate_eigen}, Proposition~\ref{Jacobian_Prop} sheds light on the algorithmic convergence rate of the DMS algorithm. As $h\to 0$ and $nh^q \to \infty$, $\nabla F(\bm{m}) \to 0$ and the convergence around any local mode $\bm{m}$ becomes superlinear. Yet, it will never be quadratic because $\nabla F(\bm{m}) \neq 0$ for any positive bandwidth $h$ and sample size $n$. In practice, one can select a small bandwidth to attain the linear convergence when the amount of available data is sufficient. Our result coincides with the findings of \cite{MS_EM2007} in Euclidean data and \cite{DirMS2020} in directional data.

\section{Discussions}
\label{Sec:Discussion}

This paper has shown that the DMS algorithm is an EM algorithm for the von Mises kernel and a generalized EM algorithm for general kernels. 
Similar to Euclidean KDE and mean shift algorithm, our results hold when each data point have a different weight and concentration parameter in the underlying directional KDE. Under the (generalized) EM representation, we have provided a new proof of the DMS algorithm's ascending property and validated the global convergence of its iterative sequences on the unit hypersphere $\Omega_{q}$. Finally, we have discussed the rate of convergence of the DMS algorithm by analyzing the Jacobian of its fixed-point iterative function.

Within a small scope, our work suggests several potential extensions or improvements of the DMS algorithm. First, instead of performing a single iteration step \eqref{EM_MS_iter} as the DMS algorithm, one may iterate \eqref{mu_fixed_point} with a fixed $\bm{\mu}^{(t)}$ for several steps before replacing $\bm{\mu}^{(t)}$ by the resulting output $\bm{\mu}^{(t+1)}$ in \eqref{mu_fixed_point}. These extra iterations in the M-step may give rise to faster convergence of the DMS algorithm. However, one needs to propose a strategy to handle zero divisions in \eqref{EM_MS_iter} when implementing this extension with truncated kernels as \eqref{convex_kernel}. Second, our explicit formula \eqref{Jacobian_F} for the Jacobian of the fixed-point iterative function can help accelerate the DMS algorithm via Aitken's or other related methods (Section 4.8 in \cite{EM_Extensions2008}). Finally, other practical guidelines for improving the regular mean shift algorithm based on its tie to the (generalized) EM algorithm \citep{MS_EM2007} may be adoptable to the DMS scenario.

More broadly, our work may have some potential impacts to the nonconvex optimization theory on manifolds. The existing convergence theory of a gradient ascent/descent type method with a nonconvex objective function (as in our DMS scenario) is often limited to some small neighborhoods of its maximizers/minimizers, especially when the objective function is supported on a compact manifold. Moreover, the ascending property of the objective function along the gradient ascent/descent path is not generally guaranteed. Our research points to a possible avenue of establishing the ascending property and global convergence theory for a broader class of gradient-based methods, that is, connecting them to the well-studied EM algorithm. Upon this connection, our technique of projecting a lower-dimensional manifold onto a higher-dimensional one could be potentially useful.

\acks{YC is supported by NSF DMS - 1810960 and DMS - 195278, NIH U01 - AG0169761.}

\bibliography{Bib_EM_MS}
	
\appendix

\section{Basins of Attraction for Directional Mean Shift Algorithm}
\label{Appendix:Basin_of_attr}

We simulate the dataset $\{\bm{X}_1,...,\bm{X}_{1000}\} \subset \Omega_2$ utilized by the DMS algorithm from a mixture of vMF densities as in \cite{DirMS2020}: 
$$f(\bm{x}) = 0.3\cdot f_{\text{vMF}}(\bm{x};\bm{\mu}_1,\kappa_1)+ 0.3\cdot f_{\text{vMF}}(\bm{x};\bm{\mu}_2,\kappa_2) + 0.4\cdot f_{\text{vMF}}(\bm{x};\bm{\mu}_3,\kappa_3)$$
with $\bm{\mu}_1 \approx (-0.35, -0.61,-0.71)$, $\bm{\mu}_2 \approx (0.5,0,0.87)$, $\bm{\mu}_3=(-0.87,0.5,0)$ (or $[-120^{\circ},-45^{\circ}]$, $[0^{\circ},60^{\circ}]$, $[150^{\circ},0^{\circ}]$ in their spherical [longitude, latitude] coordinates), and $\kappa_1=\kappa_2=8$, $\kappa_3=5$. There are three local modes of the underlying density $f$. When we apply the directional KDE \eqref{Dir_KDE} and corresponding DMS algorithm to this dataset, the bandwidth parameter $h$ is chosen using the rule-of-thumb selector (Proposition 2 in \cite{Exact_Risk_bw2013}) as $h\approx 0.384$, and the algorithms are terminated when $\max_i\norm{\bm{X}_i^{(t+1)} - \bm{X}_i^{(t)}}_2 < \epsilon=10^{-7}$. 

\begin{figure}
	\captionsetup[subfigure]{justification=centering}
	\begin{subfigure}[t]{.5\textwidth}
		\centering
		\includegraphics[width=0.9\linewidth]{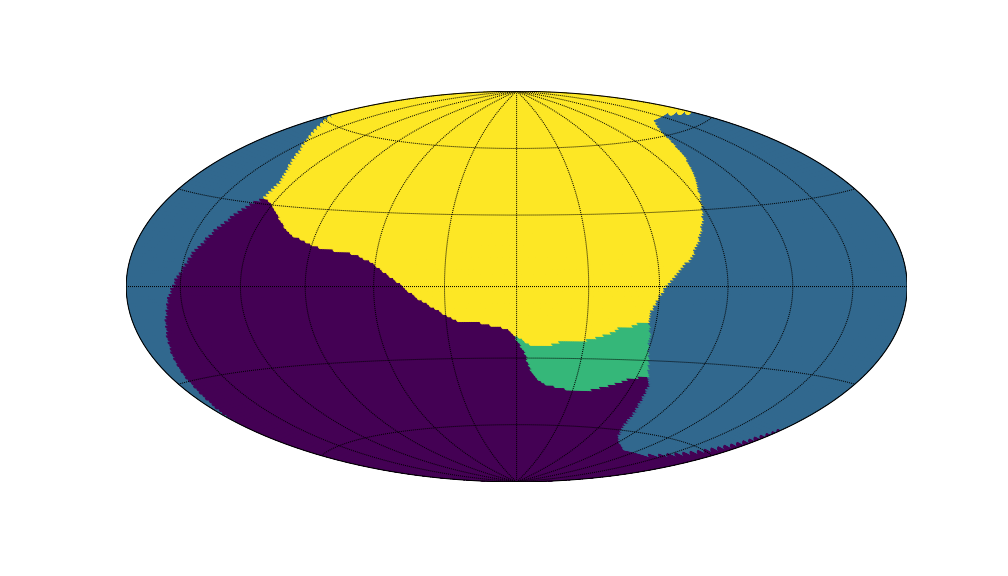}
		\caption{DMS with
			$L(r)=(1-r^2)\cdot \mathbbm{1}_{\{0\leq r \leq 1\}}$}
	\end{subfigure}%
	\hspace{1em}
	\begin{subfigure}[t]{.5\textwidth}
		\centering
		\includegraphics[width=0.9\linewidth]{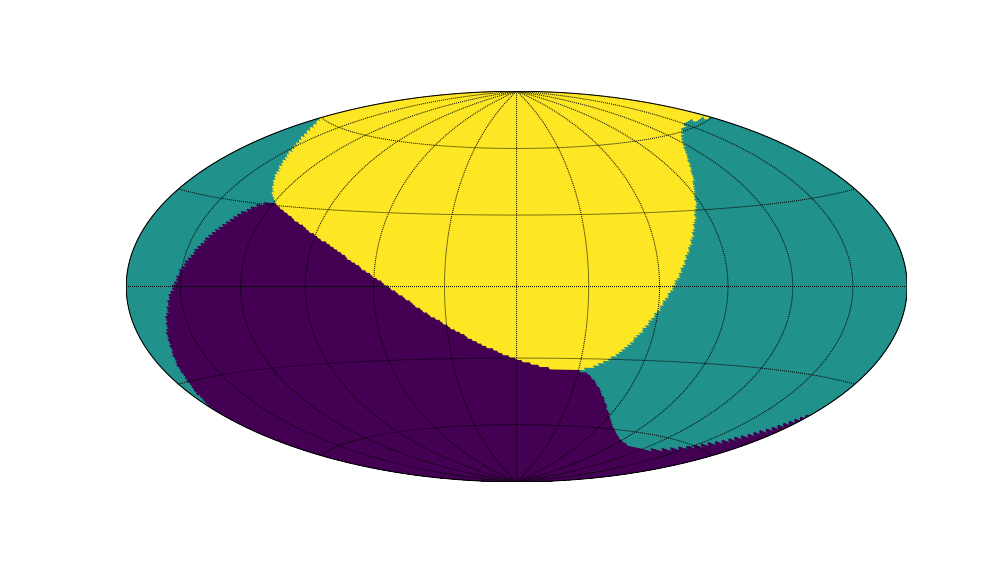}
		\caption{DMS with von Mises kernel $L(r)=e^{-r}$}
	\end{subfigure}
	\caption{Basins of attractions for DMS algorithms with two different kernels on the simulated dataset.}
	\label{fig:Basin_of_attr}
\end{figure}

Figure~\ref{fig:Basin_of_attr} displays the basins of attraction (or convergence domains) for DMS algorithms with truncated convex kernel $L(r)=(1-r^2)\cdot \mathbbm{1}_{\{0\leq r \leq 1\}}$ and von Mises kernel $L(r)=e^{-r}$. 
Each point on $\Omega_2$ except for a set with Lebesgue measure zero converges to a local mode $\hat{f}_h$ as the corresponding DMS algorithms are stopped.

\section{EM Algorithm on a Mixture of vMF Distributions}
\label{Appendix:Detail_EM_vMF}

We consider a directional dataset $\mathcal{Y}=\{\bm{Y}_1,....,\bm{Y}_N\} \subset \Omega_q$ consisting of independently and identically distributed (IID) samples from a mixture of $M$ vMF distributions, whose density is 
\begin{equation}
\label{vMF_Mixture}
f(\bm{y}|\Theta) = \sum_{j=1}^M \alpha_j \cdot f_{\text{vMF}}(\bm{y}|\bm{\mu}_j, \kappa_j).
\end{equation}
with the set of parameters as 
$$\Theta=\left\{(\alpha_j, \bm{\mu}_j,\kappa_j): \sum_{j=1}^M \alpha_j=1, \alpha_j \geq 0, \bm{\mu}_j^T\bm{\mu}_j=1,\kappa_j \geq 0, j=1,...,M \right\}.$$ 
Let $\mathcal{Z}=\{Z_1,...,Z_N\}$ be the corresponding set of hidden random variables indicating the mixture component from which the data points are sampled. In particular, $Z_i=j$ if $\bm{Y}_i$ is sampled from $f_{\text{vMF}}(\bm{y}|\bm{\mu}_j, \kappa_j)$. The EM framework for obtaining the maximum likelihood estimates for parameters in $\Theta$ given the observed data $\mathcal{Y}=\{\bm{Y}_1,....,\bm{Y}_N\} \subset \Omega_q$ is as follows \citep{spherical_EM}:\\

\noindent $\bullet$ {\bf E-Step}: Assuming that the values in $\mathcal{Z}=\{Z_1,...,Z_N\}$ are known, the complete log-likelihood is given by
\begin{equation}
\label{complete_ll}
\log P(\mathcal{Y}, \mathcal{Z}|\Theta) = \sum_{i=1}^N \sum_{j=1}^M \left[\mathbbm{1}_{\{Z_i=j\}} \log \alpha_j + \mathbbm{1}_{\{Z_i=j\}}  \log f_{\text{vMF}}(\bm{Y}_i|\bm{\mu}_j, \kappa_j) \right],
\end{equation}
where $\mathbbm{1}_A$ is the indicator function of the set $A$. For the given $(\mathcal{Y}, \Theta^{(t)})$ at the $t$-th iteration, the E-step takes the expectation of the complete log-likelihood \eqref{complete_ll} with respect to the current posterior conditional distribution of $\mathcal{Z}|\left(\mathcal{Y}, \Theta^{(t)} \right)$ as:
\begin{align}
\label{Q_function}
\begin{split}
Q(\Theta|\Theta^{(t)}) & \equiv \mathbb{E}_{\mathcal{Z}|\left(\mathcal{Y}, \Theta^{(t)} \right)}\left[\log P(\mathcal{Y},\mathcal{Z}|\Theta) \right]\\
&= \sum_{j=1}^M \sum_{i=1}^N \left(\log \alpha_j \right) \cdot P\left(Z_i=j|\mathcal{Y},\Theta^{(t)} \right) \\
&\quad + \sum_{j=1}^M \sum_{i=1}^N \left[\log f_{\text{vMF}}(\bm{Y}_i|\bm{\mu}_j,\kappa_j) \right]\cdot P\left(Z_i=j|\mathcal{Y},\Theta^{(t)} \right).\\
\end{split}
\end{align}
By Bayes' theorem, the posterior distribution $P\left(Z_i=j|\mathcal{Y},\Theta^{(t)} \right)$ is calculated as
\begin{align*}
\begin{split}
P\left(Z_i=j|\mathcal{Y},\Theta^{(t)} \right) &= \frac{P(\mathcal{Y}|Z_i=j,\Theta^{(t)}) \cdot P(Z_i=j|\Theta^{(t)})}{P(\mathcal{Y}|\Theta^{(t)})}\\
&\stackrel{\text{(**)}}{=} \frac{P(\bm{Y}_i|Z_i=j,\Theta^{(t)}) \cdot P(Z_i=j|\Theta^{(t)})}{P(\bm{Y}_i|\Theta^{(t)})} = P\left(Z_i=j|\bm{Y}_i,\Theta^{(t)} \right)\\
&= \frac{\alpha_j \cdot f_{\text{vMF}}\left(\bm{Y}_i|\bm{\mu}_j^{(t)},\kappa_j^{(t)} \right) }{\sum_{\ell=1}^M \alpha_{\ell} \cdot f_{\text{vMF}}\left(\bm{Y}_i|\bm{\mu}_{\ell}^{(t)}, \kappa_{\ell}^{(t)}\right)},
\end{split}
\end{align*} 
where we utilize the IID assumption on $\mathcal{Y}=\{\bm{Y}_1,....,\bm{Y}_N\}$ and the fact that $\bm{Y}_i$ is independent of $Z_k$ when $i\neq k$ to obtain (**).\\

\noindent $\bullet$ {\bf M-Step}: In the M-step, we solve for $\Theta^{(t+1)}$ that maximizes the function $Q(\Theta|\Theta^{(t)})$. Notice that the maximization of $\{\alpha_j\}_{j=1}^M$ and $\{\bm{\mu}_j,\kappa_j \}_{j=1}^M$ can be done separately as terms containing $\{\alpha_j\}_{j=1}^M$ and $\{\bm{\mu}_j,\kappa_j \}_{j=1}^M$ are unrelated in \eqref{Q_function}. To maximize $Q(\Theta|\Theta^{(t)})$ with respect to $\{\alpha_j\}_{j=1}^M$ under the constraint $\sum_{j=1}^M \alpha_j=1$, we introduce a Lagrangian multiplier $\lambda$ and obtain the following Lagrangian function:
$$\mathcal{L}_1(\alpha_1,...,\alpha_M,\lambda) = \sum_{j=1}^M \sum_{i=1}^N \left(\log \alpha_j \right) \cdot P\left(Z_i=j|\mathcal{Y},\Theta^{(t)} \right) + \lambda\left(\sum_{j=1}^M \alpha_j-1 \right).$$
Taking the partial derivatives of $\mathcal{L}_1$ with respect to each $\alpha_j$ and $\lambda$ and setting them to zero yield that (with some rearrangements)
$$\sum_{i=1}^N P(Z_i=j|\mathcal{Y}, \Theta^{(t)}) + \alpha_j\lambda = 0 \quad \text{ and } \quad \sum_{j=1}^M \alpha_j-1=0.$$
On summing both sides of the first equation over $j=1,...,M$ we find that $\lambda_{\text{opt}}=-N$ and hence,
\begin{equation}
\label{alpha_iter}
\alpha_j^{(t+1)} = \frac{1}{N} \sum_{i=1}^N P\left(Z_i=j|\mathcal{Y},\Theta^{(t)} \right)
\end{equation}
for $j=1,...,M$. Now, we maximize $Q(\Theta|\Theta^{(t)})$ with respect to $\{\bm{\mu}_j,\kappa_j \}_{j=1}^M$ under the constraints $\bm{\mu}_j^T \bm{\mu}_j=1$ and $\kappa_j \geq 0$ for $j=1,...,M$. Let $\lambda_j$ be the Lagrangian multiplier associated with each equality constraint $\bm{\mu}_j^T \bm{\mu}_j=1$. The Lagrangian function is given by\footnote{As mentioned in Section A.1 in \cite{spherical_EM}, we should have introduced a Lagrangian multiplier for each inequality constraint $\kappa_j\geq 0$, and work with the necessary KKT conditions. However, if $\kappa_j=0$ then $f_{\text{vMF}(\bm{y}|\bm{\mu}_j,\kappa_j)}$ is the uniform density on $\Omega_q$, and if $\kappa_j>0$ then the multiplier for the inequality constraint has to be zero by the KKT condition. Thus, the Lagrangian in \eqref{Lagrangian_2} is adequate.}

\begin{align}
\label{Lagrangian_2}
\begin{split}
\mathcal{L}_2\left(\{\bm{\mu}_j,\kappa_j,\lambda_j\}_{j=1}^M \right) &= \sum_{j=1}^M \sum_{i=1}^N \left[\log f_{\text{vMF}}(\bm{Y}_i|\bm{\mu}_j,\kappa_j) \right]\cdot P\left(Z_i=j|\mathcal{Y},\Theta^{(t)} \right) + \sum_{j=1}^M \lambda_j (1-\bm{\mu}_j^T\bm{\mu}_j)\\
&= \sum_{j=1}^M \left\{\sum_{i=1}^N \left[\log C_q(\kappa_j) + \kappa_j\bm{\mu}_j^T \bm{Y}_i \right]\cdot  P\left(Z_i=j|\mathcal{Y},\Theta^{(t)} \right) +\lambda_j (1-\bm{\mu}_j^T\bm{\mu}_j) \right\}.
\end{split}
\end{align}
Taking partial derivatives of \eqref{Lagrangian_2} with respect to $\{\bm{\mu}_j,\kappa_j,\lambda_j\}_{j=1}^M$ and setting them to zero, we have that
\begin{align}
\label{Lagrangian_2_Deri1}
\kappa_j \sum_{i=1}^N \bm{Y}_i \cdot P(Z_i=j|\mathcal{Y},\Theta^{(t)}) -2\lambda_j \bm{\mu}_j &=0,\\
\label{Lagrangian_2_Deri2}
\frac{C_q'(\kappa_j)}{C_q(\kappa_j)} \sum_{i=1}^N P(Z_i=j|\mathcal{Y},\Theta^{(t)})+ \bm{\mu}_j \sum_{i=1}^N \bm{Y}_i \cdot P(Z_i=j|\mathcal{Y}, \Theta^{(t)}) &= 0,\\
\label{Lagrangian_2_Deri3}
1-\bm{\mu}_j^T \bm{\mu}_j &= 0,
\end{align}
for $j=1,...,M$. Combining \eqref{Lagrangian_2_Deri1} and \eqref{Lagrangian_2_Deri3} shows that for $j=1,...,M$,
$$\hat{\lambda}_j = \frac{\kappa_j}{2} \norm{\sum_{i=1}^N \bm{Y}_i \cdot P\left(Z_i=j|\mathcal{Y}, \Theta^{(t)}\right)}_2,$$
\begin{equation}
\label{mu_iter}
\bm{\mu}_j^{(t+1)} = \frac{\sum_{i=1}^N \bm{Y}_i \cdot P(Z_i=j|\mathcal{Y}, \Theta^{(t)})}{\norm{\sum_{i=1}^N \bm{Y}_i \cdot P(Z_i=j|\mathcal{Y}, \Theta^{(t)})}_2}.
\end{equation}
Substituting \eqref{mu_iter} into \eqref{Lagrangian_2_Deri2} gives us that
\begin{equation}
\label{kappa_iter}
-\frac{C_q'\left(\kappa_j^{(t+1)} \right)}{C_q\left(\kappa_j^{(t+1)} \right)} = \frac{\mathcal{I}_{\frac{q+1}{2}}\left(\kappa_j^{(t+1)} \right)}{\mathcal{I}_{\frac{q-1}{2}}\left(\kappa_j^{(t+1)} \right)} = \frac{\norm{\sum_{i=1}^N \bm{Y}_i \cdot P(Z_i=j|\mathcal{Y}, \Theta^{(t)})}_2}{\sum_{i=1}^N P(Z_i=j|\mathcal{Y},\Theta^{(t)})}.
\end{equation}
The first equality in \eqref{kappa_iter} follows from the below calculations (recall \eqref{vMF_density}):
\begin{align*}
C_q'(\kappa) &= \frac{\left(\frac{q-1}{2} \right) \kappa^{\frac{q-3}{2}}(2\pi)^{\frac{q+1}{2}} \mathcal{I}_{\frac{q-1}{2}}(\kappa) - (2\pi)^{\frac{q+1}{2}} \mathcal{I}_{\frac{q-1}{2}}'(\kappa) \kappa^{\frac{q-1}{2}}}{(2\pi)^{q+1}\left[\mathcal{I}_{\frac{q-1}{2}}(\kappa) \right]^2}\\
&= \frac{\kappa^{\frac{q-1}{2}}}{(2\pi)^{\frac{q+1}{2}} \mathcal{I}_{\frac{q-1}{2}}(\kappa)}\left[\frac{q-1}{2\kappa} - \frac{\mathcal{I}_{\frac{q-1}{2}}'(\kappa)}{\mathcal{I}_{\frac{q-1}{2}}(\kappa)} \right]
\end{align*}
and thus,
\begin{align*}
\frac{C_q'(\kappa)}{C_q(\kappa)} &= \frac{q-1}{2\kappa} - \frac{\mathcal{I}_{\frac{q-1}{2}}'(\kappa)}{\mathcal{I}_{\frac{q-1}{2}}(\kappa)} = \frac{q-1}{2\kappa} - \frac{\left(\frac{q-1}{2\kappa} \right) \mathcal{I}_{\frac{q-1}{2}}(\kappa) + \mathcal{I}_{\frac{q+1}{2}}(\kappa)}{\mathcal{I}_{\frac{q-1}{2}}(\kappa)} =- \frac{\mathcal{I}_{\frac{q+1}{2}}(\kappa)}{\mathcal{I}_{\frac{q-1}{2}}(\kappa)},
\end{align*}
where we use the well-know recurrence relation $\kappa\mathcal{I}_s'(\kappa) = s \mathcal{I}_s(\kappa) + \kappa \mathcal{I}_{s+1}(\kappa)$; see, for instance, Section 9.6.26 in \cite{Handbook_math1974}.\\
The updating formula \eqref{kappa_iter} involves the modified Bessel function of the first kind and has no closed-form solutions for $\kappa_j^{(t+1)}$. However, \cite{spherical_EM,Sra2011} propose an approximation formula for $\kappa$ in terms of $A_q(\kappa) \equiv \frac{\mathcal{I}_{\frac{q+1}{2}}(\kappa)}{\mathcal{I}_{\frac{q-1}{2}}(\kappa)}$ as
\begin{equation}
\label{kappa_approx}
\hat{\kappa} \approx \frac{(q+1)A_q(\kappa) - A_q(\kappa)^3}{1- A_q(\kappa)^2}.
\end{equation}
Therefore, we may update $\kappa_j^{(t+1)}$ in the E-step by
\begin{equation}
\label{kappa_iter_approx}
\kappa_j^{(t+1)} = \frac{(q+1)A_q(\kappa^{(t+1)}) - A_q(\kappa^{(t+1)})^3}{1- A_q(\kappa^{(t+1)})^2}
\end{equation}
with $A_q(\kappa^{(t+1)}) = \frac{\norm{\sum_{i=1}^N \bm{Y}_i \cdot P(Z_i=j|\mathcal{Y}, \Theta^{(t)})}_2}{\sum_{i=1}^N P(Z_i=j|\mathcal{Y},\Theta^{(t)})}$.\\

The above EM algorithm involves only soft assignments for hidden variables $\mathcal{Z}$.
The analysis of hard assignments for hidden variables and the subsequent discussion on their connections with the spherical $k$-Means \citep{Spherical_Kmeans2001} clustering can be found in \cite{spherical_EM}.

\section{Naive Method: Introducing a Displacement Parameter}
\label{Appendix:EM_displacement}

We show in this section that introducing a displacement parameter to define an artificial maximum likelihood as in \cite{MS_EM2007} is (probably) not a promising direction under our directional data scenario. For simplicity, we only consider the von Mises kernel case. Given a directional dataset $\{\bm{X}_1,..,\bm{X}_n\} \subset \Omega_q$, we may define the mixture model by subtracting a displacement parameter $\bm{\mu}\in \Omega_q$ from each mean $\{\bm{X}_1,..,\bm{X}_n\}$ of the mixture component as:
\begin{equation}
\label{mixture_den_dis}
f_{\text{MOvMF}}(\bm{y}) = \sum_{i=1}^n \alpha_i \cdot C_q(\kappa_i,\bm{X}_i-\bm{\mu}) \cdot \exp\left[\kappa_i(\bm{X}_i-\bm{\mu})^T\bm{y}\right].
\end{equation}
The key difference between the mixture density here and the one in \eqref{Mixture_Density} is that the normalizing constant in each component of \eqref{mixture_den_dis} depends not only on the (fixed) concentration parameter $\kappa_i$ but also our newly introduced displacement parameter $\bm{\mu}$. The reason is that a unit vector subtracted by another unit vector is no longer a unit vector on $\Omega_q$. Hence,
\begin{align*}
C_q(\kappa,\bm{X}-\bm{\mu})^{-1} &= \int_{\Omega_q} \exp\left[\kappa (\bm{X}-\bm{\mu})^T \bm{y} \right] \omega_q(d\bm{y})\\ 
&= \omega_{q-1}(\Omega_{q-1}) \int_{-1}^1 \exp\left(\kappa t \norm{\bm{X}-\bm{\mu}}_2 \right) \cdot (1-t^2)^{\frac{q}{2}-1} dt\\
&= \frac{(2\pi)^{\frac{q+1}{2}}}{\left(\kappa \norm{\bm{X}-\bm{\mu}}_2\right)^{\frac{q-1}{2}}} \cdot \mathcal{I}_{\frac{q-1}{2}}\left( \kappa \norm{\bm{X}-\bm{\mu}}_2\right)
\end{align*}
This extra dependence becomes a main obstacle to connecting the corresponding EM algorithm with the directional mean shift algorithm, because one can easily verify that the derivative of the resulting Q-function with respect to $\bm{\mu}$ in this case would involve the derivative of modified Bessel functions of the first kind. There is no closed-form expression for such derivative in general.

\section{Proofs of Theorems and Propositions}

\subsection{Proof of Theorem~\ref{Q_fun_increase}}
\label{Appendix:Q_fun_increase_pf}

\begin{customthm}{1}
	Under condition (C2), the Q-function \eqref{Q_function_MS2} satisfies $Q(\bm{\mu}^{(t+1)}|\bm{\mu}^{(t)}) \geq Q(\bm{\mu}^{(t)}|\bm{\mu}^{(t)})$ 
	where $\bm{\mu}^{(t+1)} = -\frac{\sum_{i=1}^n \bm{X}_i L'\left(\frac{1-\bm{X}_i^T\bm{\mu}^{(t)}}{h^2}\right)}{\norm{\sum_{i=1}^n \bm{X}_i L'\left(\frac{1-\bm{X}_i^T\bm{\mu}^{(t)}}{h^2}\right)}_2}$ is from the iteration \eqref{EM_MS_iter}.
\end{customthm}

\begin{proof}
	As we take $\kappa_i=\frac{1}{h^2}$ and $\alpha_i=\frac{1}{n}$ for all $i=1,...,n$ in equations \eqref{Q_function_MS2} and \eqref{post_hidden_MS}, the Q-function reduces to
	\begin{equation}
	\label{Q_function_MS_nh}
	Q(\bm{\mu}|\bm{\mu}^{(t)}) = \sum_{i=1}^n \left[\log\frac{C_{1/h^2,q+1,L}}{n} + \log L\left(\frac{1-\bm{X}_i^T\bm{\mu}}{h^2} \right) \right] \cdot P(Z_1=i|\bm{Y}_1,\bm{\mu}^{(t)}),
	\end{equation}
	and the posterior $Z_1|\bm{Y}_1,\bm{\mu}^{(t)}$ becomes
	\begin{equation}
	\label{post_hidden_MS_nh}
	P(Z_1=i|\bm{Y}_1,\bm{\mu}^{(t)}) = \frac{L\left(\frac{1-\bm{X}_i^T\bm{\mu}^{(t)}}{h^2} \right)}{\sum_{\ell=1}^n L\left(\frac{1-\bm{X}_{\ell}^T\bm{\mu}^{(t)}}{h^2} \right)}.
	\end{equation}
	The differentiability and convexity of $L$ under (C2) indicates that
	\begin{equation}
	\label{convex_inq}
	L(y)-L(x) \geq L'(x) \cdot (y-x),
	\end{equation}
	where $L'(x)$ is replaced by any subgradient at non-differentiable points. As $0<L(0)<\infty$, we have that
	\begin{align*}
	&Q(\bm{\mu}^{(t+1)}|\bm{\mu}^{(t)}) - Q(\bm{\mu}^{(t)}|\bm{\mu}^{(t)})\\ 
	&= \sum_{i=1}^n \left[\log L\left(\frac{1-\bm{X}_i^T\bm{\mu}^{(t+1)}}{h^2} \right) - \log L\left(\frac{1-\bm{X}_i^T\bm{\mu}^{(t)}}{h^2} \right) \right] \cdot \frac{L\left(\frac{1-\bm{X}_i^T\bm{\mu}^{(t)}}{h^2} \right)}{\sum_{\ell=1}^n L\left(\frac{1-\bm{X}_{\ell}^T\bm{\mu}^{(t)}}{h^2} \right)}\\
	&\stackrel{\text{(i)}}{\geq} \sum_{i=1}^n \frac{1}{L\left(\frac{1-\bm{X}_i^T\bm{\mu}^{(t+1)}}{h^2} \right)} \cdot \left[ L\left(\frac{1-\bm{X}_i^T\bm{\mu}^{(t+1)}}{h^2} \right) - L\left(\frac{1-\bm{X}_i^T\bm{\mu}^{(t)}}{h^2} \right) \right] \cdot \frac{L\left(\frac{1-\bm{X}_i^T\bm{\mu}^{(t)}}{h^2} \right)}{\sum_{\ell=1}^n L\left(\frac{1-\bm{X}_{\ell}^T\bm{\mu}^{(t)}}{h^2} \right)}\\
	&\stackrel{\text{(ii)}}{\geq} \sum_{i=1}^n \frac{L'\left(\frac{1-\bm{X}_i^T\bm{\mu}^{(t)}}{h^2} \right)}{L\left(\frac{1-\bm{X}_i^T\bm{\mu}^{(t+1)}}{h^2} \right)} \cdot \left[\frac{\bm{X}_i^T\left(\bm{\mu}^{(t)} -\bm{\mu}^{(t+1)} \right)}{h^2} \right] \cdot \frac{L\left(\frac{1-\bm{X}_i^T\bm{\mu}^{(t)}}{h^2} \right)}{\sum_{\ell=1}^n L\left(\frac{1-\bm{X}_{\ell}^T\bm{\mu}^{(t)}}{h^2} \right)}\\
	&\stackrel{\text{(iii)}}{\geq} \frac{L\left(\frac{2}{h^2} \right)}{nL(0)^2}\sum_{i=1}^n L'\left(\frac{1-\bm{X}_i^T\bm{\mu}^{(t)}}{h^2} \right) \bm{X}_i^T \left(\bm{\mu}^{(t)} -\bm{\mu}^{(t+1)} \right)\\
	&\stackrel{\text{(iv)}}{=} \frac{L\left(\frac{2}{h^2} \right)}{nL(0)^2} \norm{\sum_{i=1}^n \bm{X}_i L'\left(\frac{1-\bm{X}_i^T\bm{\mu}^{(t)}}{h^2} \right)} \cdot \left(\bm{\mu}^{(t+1)}\right)^T \left(\bm{\mu}^{(t+1)} -\bm{\mu}^{(t)} \right)\\
	&\stackrel{\text{(v)}}{=} \frac{L\left(\frac{2}{h^2} \right)}{2nL(0)^2} \norm{\sum_{i=1}^n \bm{X}_i L'\left(\frac{1-\bm{X}_i^T\bm{\mu}^{(t)}}{h^2} \right)} \cdot \norm{\bm{\mu}^{(t+1)} -\bm{\mu}^{(t)}}_2 \\
	&\geq 0,
	\end{align*}
	where we use the fact that $\log y -\log x \geq \frac{y-x}{y}$ by the concavity of $\log(\cdot)$ in (i), leverage the inequality \eqref{convex_inq} in (ii), apply the monotonicity of $L$ in (iii), plug in \eqref{EM_MS_iter} in (iv), and notice that $\left(\bm{\mu}^{(t+1)} \right)^T \left(\bm{\mu}^{(t+1)} -\bm{\mu}^{(t)} \right)= \frac{1}{2} \norm{\bm{\mu}^{(t+1)} -\bm{\mu}^{(t)}}_2$ on $\Omega_q$ in (v). The result follows.
\end{proof}

\subsection{Proof of Proposition~\ref{obs_ll_Dir_KDE}}
\label{Appendix:obs_ll_Dir_KDE}

\begin{customprop}{2}
	The observed log-likelihood $\log P(\bm{Y}_1|\bm{\mu})$ is given by 
	\begin{equation}
	\label{obs_log_likelihood2}
	\log P(\bm{Y}_1|\bm{\mu}) = \log\left[\sum_{\ell=1}^n \alpha_{\ell} \cdot C_{\kappa_{\ell},q+1,L} \cdot L\left[\kappa_{\ell} \left(1-\bm{X}_{\ell}^T\bm{\mu} \right)\right] \right].
	\end{equation}
	In particular, when $\alpha_{\ell}=\frac{1}{n}$ and $\kappa_{\ell}=\frac{1}{h^2}$ for all $\ell=1,...,n$ (i.e., simplifying to the DMS case), it implies that
	\begin{equation}
	\label{obs_ll_log_density3}
	\log P(\bm{Y}_1|\bm{\mu}) = \log \hat{f}_h(\bm{\mu}) + \log \frac{C_{1/h^2,q+1,L}}{c_{h,q,L}},
	\end{equation}
	where $c_{h,q,L}$ is defined in \eqref{asym_norm_const} and $C_{\kappa_i,q+1,L}$ is elucidated in \eqref{norm_const}.
\end{customprop}

\begin{proof}
	Notice that the observed log-likelihood in the EM framework of Section \ref{Sec:Detailed_deri}, after we take $\bm{Y}_1=(0,...,0,1)^T \in \Omega_{q+1}$, is given by
	\begin{align}
	\label{obs_ll}
	\begin{split}
	\log P(\bm{Y}_1|\bm{\mu}) &= \log P(\bm{Y}_1,Z_1|\bm{\mu}) - \log P(Z_1|\bm{Y}_1,\bm{\mu})\\
	&= \sum_{i=1}^n \mathbbm{1}_{\{Z_1=i\}} \cdot \left\{\log \alpha_i + \log C_{\kappa_i,q+1,L} + \log L\left[\kappa_i(1-\bm{X}_i^T\bm{\mu})\right] \right\}\\
	& \quad - \sum_{i=1}^n \mathbbm{1}_{\{Z_1=i \}} \cdot \log P(Z_1=i|\bm{Y}_1,\bm{\mu}).
	\end{split}
	\end{align}
	Taking the expectation over $Z_1|\left(\bm{Y}_1,\bm{\mu}^{(t)}\right)$ on both sides of \eqref{obs_ll} yields that
	\begin{align}
	\label{obs_ll2}
	\begin{split}
	\log P(\bm{Y}_1|\bm{\mu}) &= \sum_{i=1}^n P\left(Z_1=i|\bm{Y}_1,\bm{\mu}^{(t)}\right) \cdot \left\{\log \alpha_i + \log C_{\kappa_i,q+1,L} + \log L\left[\kappa_i(1-\bm{X}_i^T\bm{\mu})\right] \right\}\\
	& \quad - \sum_{i=1}^n P\left(Z_1=i|\bm{Y}_1,\bm{\mu}^{(t)} \right) \cdot \log P(Z_1=i|\bm{Y}_1,\bm{\mu}).
	\end{split}
	\end{align}
	After plugging \eqref{post_hidden_MS} into \eqref{obs_ll2}, we obtain that
	\begin{align*}
	&\log P(\bm{Y}_1|\bm{\mu})\\
	&= \sum_{i=1}^n P\left(Z_1=i|\bm{Y}_1,\bm{\mu}^{(t)} \right) \cdot \Bigg[\log \alpha_i + \log C_{\kappa_i,q+1,L} + \log L\left[\kappa_i(1-\bm{X}_i^T\bm{\mu}) \right] \\
	& \quad\quad \quad- \log \frac{\alpha_i \cdot C_{\kappa_i,q+1,L} \cdot L\left[\kappa_i(1-\bm{X}_i^T \bm{\mu}) \right]}{\sum_{\ell=1}^n \alpha_{\ell} \cdot C_{\kappa_{\ell},q+1,L} \cdot L\left[\kappa_{\ell} \left(1-\bm{X}_{\ell}^T\bm{\mu} \right)\right]} \Bigg]\\
	&= \sum_{i=1}^n P\left(Z_1=i|\bm{Y}_1,\bm{\mu}^{(t)} \right) \cdot \log \left[\sum_{\ell=1}^n \alpha_{\ell} \cdot C_{\kappa_{\ell},q+1,L} \cdot L\left[\kappa_{\ell} \left(1-\bm{X}_{\ell}^T\bm{\mu} \right)\right] \right]\\
	&= \log\left[\sum_{\ell=1}^n \alpha_{\ell} \cdot C_{\kappa_{\ell},q+1,L} \cdot L\left[\kappa_{\ell} \left(1-\bm{X}_{\ell}^T\bm{\mu} \right)\right] \right].
	\end{align*}
	When $\alpha_{\ell}=\frac{1}{n}$ and $\kappa_{\ell}=\frac{1}{h^2}$ for all $\ell =1,...,n$ (i.e., simplifying to the DMS form), the above equation indicates that
	\begin{align}
	\label{obs_ll_log_density2}
	\begin{split}
	\log P(\bm{Y}_1|\bm{\mu}) &= \log \left[\frac{c_{h,q,L}}{n} \sum_{i=1}^n L\left(\frac{1-\bm{X}_i^T\bm{\mu}}{h^2} \right) \right] + \log \frac{C_{1/h^2,q+1,L}}{c_{h,q,L}} \\
	&= \log \hat{f}_h(\bm{\mu}) + \log \frac{C_{1/h^2,q+1,L}}{c_{h,q,L}},
	\end{split}
	\end{align}
	where $c_{h,q,L}$ is defined in \eqref{asym_norm_const} and $C_{\kappa_i,q+1,L}$ is elucidated in \eqref{norm_const}. The results follow.
\end{proof}

\subsection{Proof of Theorem~\ref{EM_Ascending}}
\label{Appendix:EM_Ascending}

\begin{customthm}{3}
	Let $\left\{\bm{\mu}^{(t)} \right\}_{t=0}^{\infty}$ be the path of successive points defined by the DMS algorithm \eqref{Dir_MS}. Under condition (C2), the sequence $\left\{\hat{f}_h(\bm{\mu}^{(t)}) \right\}_{t=0}^{\infty}$ is non-decreasing and thus converges.
\end{customthm}

\begin{proof}
	The difference of observed log-likelihoods \eqref{obs_ll2} between two consecutive steps of our DMS (or generalized EM) iteration is calculated as
	\begin{align}
	\label{obs_ll_diff}
	\begin{split}
	&\log P\left(\bm{Y}_1|\bm{\mu}^{(t+1)} \right) - \log P\left(\bm{Y}_1|\bm{\mu}^{(t)} \right) \\
	&\stackrel{\text{(i)}}{=} Q\left(\bm{\mu}^{(t+1)}|\bm{\mu}^{(t)} \right) - \sum_{i=1}^n P\left(Z_1=i|\bm{Y}_1,\bm{\mu}^{(t)} \right) \cdot \log P\left(Z_1=i|\bm{Y}_1,\bm{\mu}^{(t+1)} \right)\\
	&\quad - Q\left(\bm{\mu}^{(t)}|\bm{\mu}^{(t)} \right) + \sum_{i=1}^n P\left(Z_1=i|\bm{Y}_1,\bm{\mu}^{(t)} \right) \cdot \log P\left(Z_1=i|\bm{Y}_1,\bm{\mu}^{(t)} \right)\\
	&= Q\left(\bm{\mu}^{(t+1)}|\bm{\mu}^{(t)} \right) - Q\left(\bm{\mu}^{(t)}|\bm{\mu}^{(t)} \right) - \sum_{i=1}^n P\left(Z_1=i|\bm{Y}_1,\bm{\mu}^{(t)} \right) \cdot \log \frac{P\left(Z_1=i|\bm{Y}_1,\bm{\mu}^{(t+1)} \right)}{P\left(Z_1=i|\bm{Y}_1,\bm{\mu}^{(t)} \right)}\\
	&\stackrel{\text{(ii)}}{\geq} Q\left(\bm{\mu}^{(t+1)}|\bm{\mu}^{(t)} \right) - Q\left(\bm{\mu}^{(t)}|\bm{\mu}^{(t)} \right) - \underbrace{\log\left[\sum_{i=1}^n P\left(Z_1=i|\bm{Y}_1,\bm{\mu}^{(t+1)} \right) \right]}_{=0}\\
	& \stackrel{\text{(iii)}}{\geq} 0,
	\end{split}
	\end{align}
	where we recall \eqref{Q_function_MS2} in (i), apply Jensen's inequality in (ii), and utilize Theorem~\ref{Q_fun_increase} in (iii). The inequality (ii) is strict except when $\bm{\mu}^{(t+1)} = \bm{\mu}^{(t)}$.
	Therefore, combining \eqref{obs_ll_log_density2} with \eqref{obs_ll_diff} show that
	$$\log \hat{f}_h(\bm{\mu}^{(t+1)} ) - \log\hat{f}_h(\bm{\mu}^{(t)}) = \log P\left(\bm{Y}_1|\bm{\mu}^{(t+1)} \right) - \log P\left(\bm{Y}_1|\bm{\mu}^{(t)} \right) \geq 0.$$
	As $\hat{f}_h$ is differentiable on the compact set $\Omega_q$ under (C2), $\left\{\hat{f}_h(\bm{\mu}^{(t)}) \right\}_{t=0}^{\infty}$ is thus bounded from above and converges.
\end{proof}

\subsection{Conditions on the Convergence of Generalized EM Sequences and Proof of Theorem~\ref{EM_conv_new}}
\label{Appendix:EM_conv}

For the convergence of (generalized) EM sequences to stationary points (in most cases, local maxima), \cite{EM_Jeff1983} made the following assumptions on the parameter space $\Psi$ and observed log-likelihood function $\ell(\phi)$:
\begin{itemize}
	\item {\bf (A1)} $\Psi$ is a subset of the ambient Euclidean space,
	\item {\bf (A2)} $\Psi_0=\left\{\phi \in \Psi: \ell(\phi) \geq \ell(\phi_0) \right\}$ is compact for any $\ell(\phi_0) > -\infty$, where $\phi_0$ is the initial parameter for the (generalized) EM algorithm,
	\item {\bf (A3)} $\ell$ is continuous in $\Psi$ and differentiable in the interior of $\Psi$.
\end{itemize}
Let $\phi \to M(\phi)$ be the point-to-set map defined by a generalized EM iteration such that the Q-function\footnote{The Q-function is defined as the expectation of the complete log-likelihood over the posterior distribution of hidden variables given the values of observed variables and current parameters.} satisfies
$$Q(\phi'|\phi) \geq Q(\phi|\phi) \quad \text{ for any } \phi'\in \Psi.$$
Denote by $\mathcal{S}$ the set of stationary points of the observed log-likelihood $\ell$. As a consequence of the above assumptions, all the limit points of the generalized EM sequence are stationary points of the observed log-likelihood. 

\begin{lemma}[Theorem 1 of \cite{EM_Jeff1983}]
	\label{EM_conv}
	Let $\{\phi_p\}$ be a generalized EM sequence generated by $\phi_{p+1} \in M(\phi_p)$, and suppose that (i) $M$ is a closed point-to-set map over the complement of $\mathcal{S}$; (ii) $\ell(\phi_{p+1}) > \ell(\phi_p)$ for all $\phi_p \notin \mathcal{S}$. Then all the limit points of $\{\phi_p\}$ are stationary points of $\ell$, and $\ell(\phi_p)$ converges monotonically to $\ell^*=\ell(\phi^*)$ for some $\phi^* \in \mathcal{S}$.
\end{lemma}

The generalized EM iteration is usually terminated when the (Euclidean) distance between two consecutive iteration points is below some small threshold value. However, this property does not imply the convergence of the generalized EM sequence to a single point in $\mathcal{S}$; see some related discussions in \cite{EM_Jeff1983} and a counterexample in \cite{MS2007_pf}. To guarantee the convergence of a generalized EM sequence to a single point in $\mathcal{S}$ (in most cases, a local maximum), we require the discreteness of the set of local modes:
\begin{itemize}
	\item {\bf (A4)} The set of local modes of the observed log-likelihood $\ell(\phi)$ is discrete, and the modes are isolated from other critical points.
\end{itemize}

\subsubsection{Proof of Theorem~\ref{EM_conv_new}}

\begin{customthm}{4}
	Let $\left\{\bm{\mu}^{(t)} \right\}_{t=0}^{\infty}$ be the path of successive points defined by the DMS algorithm \eqref{Dir_MS}. Under conditions (C1) and (C2), the sequence $\left\{\bm{\mu}^{(t)} \right\}_{t=0}^{\infty}$ converges to a local mode of $\hat f_h$ from almost all starting point $\bm{\mu}^{(0)}\in \Omega_{q}$.
\end{customthm}

\begin{proof}
	Now, we verify in detail that the above assumptions (A1-4) for the convergence of generalized EM sequences can be satisfied by our conditions (C1) and (C2) on directional KDE $\hat{f}_h$ and kernel function $L$. Notice that the estimated parameter $\bm{\mu}$ lies on $\Omega_q$, a compact manifold in the ambient Euclidean space $\mathbb{R}^{q+1}$. By the equivalence of the observed log-likelihood and logarithm of the directional KDE (c.f., equation \eqref{obs_ll_log_density}), the constrained parameter space $\Psi_0= \left\{\bm{\mu}\in \Omega_q: \log \hat{f}_h(\bm{\mu}) \geq \log \hat{f}_h(\bm{\mu}^{(0)}) \right\}$ is also compact for any $\log \hat{f}_h(\bm{\mu}^{(0)}) > -\infty$. Moreover, under condition (C2), $\hat{f}_h$ is differentiable on $\Omega_q$. Thus, assumptions (A1-3) are naturally justified in our DMS context. Finally, it is obvious that our condition (C1) implies assumption (A4), and as long as we do not initialize the DMS algorithm within the set of local minima and saddle points, the algorithm will converge to a local mode.
\end{proof}

\subsection{Proof of Proposition~\ref{Jacobian_Prop}}
\label{Appendix:Jacobian_Prop_pf}

\begin{customprop}{5}
	For any $\bm{x}\in \Omega_q \subset \mathbb{R}^{q+1}$,
	\begin{enumerate}[label=(\alph*)]
		\item the Jacobian $\nabla F(\bm{x})$ is a symmetric matrix in $\mathbb{R}^{(q+1)\times (q+1)}$, and has its explicit form as
		\begin{equation}
		\label{Jacobian_F1}
		\nabla F(\bm{x}) = \frac{\nabla\nabla \hat{f}_h(\bm{x})}{\norm{\nabla \hat{f}_h(\bm{x})}_2} - \frac{\nabla \hat{f}_h(\bm{x}) \nabla \hat{f}_h(\bm{x})^T \nabla\nabla \hat{f}_h(\bm{x})}{\norm{\nabla \hat{f}_h(\bm{x})}_2^3}.
		\end{equation}
		In particular, $\nabla F(\bm{m}) = \frac{(\bm{I}_{q+1}-\bm{m}\bm{m}^T) \nabla \nabla \hat{f}_h(\bm{m})}{\norm{\nabla \hat{f}_h(\bm{m})}_2}$ at any local mode $\bm{m}$ of $\hat{f}_h$.
		\item we have that $\max_{i=1,...,q+1}|\Lambda_i| \to 0$ as $h\to 0$ and $nh^q\to \infty$, where $\Lambda_1,...,\Lambda_{q+1}$ are the eigenvalues of $\nabla F(\bm{m})$.
	\end{enumerate}
\end{customprop}

\begin{proof}
	(a) We leverage the form $F(\bm{x}) = \frac{\nabla\hat{f}_h(\bm{x})}{\norm{\nabla\hat{f}_h(\bm{x})}_2}$ in \eqref{Dir_MS_fix} to directly compute the explicit form of the Jacobian $F(\bm{x})$ as
	$$\nabla F(\bm{x}) = \frac{\nabla\nabla \hat{f}_h(\bm{x})}{\norm{\nabla \hat{f}_h(\bm{x})}_2} - \frac{\nabla \hat{f}_h(\bm{x}) \nabla \hat{f}_h(\bm{x})^T \nabla\nabla \hat{f}_h(\bm{x})}{\norm{\nabla \hat{f}_h(\bm{x})}_2^3},$$
	which is a symmetric matrix. At local mode $\bm{m}$, we have that $\bm{m} = \frac{\nabla \hat{f}_h(\bm{m})}{\norm{\nabla \hat{f}_h(\bm{m})}_2}$, and the Jacobian $\nabla F(\bm{m})$ is simplified as
	\begin{equation*}
	\nabla F(\bm{m}) = \frac{(\bm{I}_{q+1}-\bm{m}\bm{m}^T) \nabla \nabla \hat{f}_h(\bm{m})}{\norm{\nabla \hat{f}_h(\bm{m})}_2}, 
	\end{equation*}
	where $I_{q+1}\in \mathbb{R}^{(q+1)\times (q+1)}$ is the identity matrix. \\
	
	\noindent (b) By the definition of local modes of $\hat{f}_h$ on $\Omega_q$, the total gradient $\nabla \hat{f}_h(\bm{m})$ has zero projection onto the tangent space of $\Omega_q$ at any local mode $\bm{m}$. Thus, $\nabla \hat{f}_h(\bm{m})$ is completely determined by its radial component along $\bm{m}$. By Lemma 10 in \cite{DirMS2020}, it is obvious that $\norm{\nabla \hat{f}_h(\bm{m})}_2 \to \infty$ as $h\to 0$ and $nh^q \to \infty$. Finally, each eigenvalue of $\nabla F(\bm{m})$ will have its absolute value tending to 0 as $\norm{\nabla \hat{f}_h(\bm{m})}_2 \to \infty$. The result follows.
\end{proof}
	
	
\end{document}